
%
%
%
%
%
%
\input amstex
\documentstyle {amsppt}
\tolerance=10000
\magnification=1200
%
%
\hsize=6.7truein 
%
%
\addto\tenpoint{\normalbaselineskip=18pt\normalbaselines}
\addto\eightpoint{\normalbaselineskip=15pt\normalbaselines}
\refstyle{C}
\NoRunningHeads
\def\R{\Bbb R}
\def\Z{\Bbb Z}
\def\N{\Bbb N}
\def\T{\Bbb T}
\def\C{\Bbb C}

\def\s{\sigma}
\def\GA{G\rtimes A}
\def\GB{G\rtimes B}
\def\schwmh{\Cal S_{H}^{\sigma}(M)}
\def\Kinf{{\Cal K}^{\infty}}
\def\K{\Cal K}
\def\H{\Cal H}
\def\schwmns{{\Cal S} (M)}

\def\pa{\parallel}

\def\tensp{\hat \otimes_{\pi}}
\def\schwr{\Cal S (\R)}
\def\Atil{\tilde {A}}
\def\Botil{\tilde B}
\def\schwmg{{\Cal S}^{\s}_{G} (M)}
\topmatter
\title Spectral Invariance of Dense Subalgebras of Operator Algebras
\endtitle
\author Larry B. Schweitzer \endauthor
\address Department of Mathematics and Statistics, University of Victoria,
Victoria, British Columbia, Canada V8W 3P4 \endaddress
\email lschweit\@ sol.uvic.ca \endemail
\endtopmatter
\document
\heading Abstract \endheading
\par
We define the notion of strong spectral invariance
for a dense Fr\'echet subalgebra $A$ of a Banach algebra $B$.  We show
that if $A$ is strongly spectral invariant in a C*-algebra $B$, and $G$
is a compactly generated polynomial growth Type R Lie group,
not necessarily connected, then the
smooth crossed product $G\rtimes A$ is spectral invariant in
the C*-crossed product $G \rtimes B$.   Examples
of such groups are given by finitely generated polynomial growth
discrete groups, compact or connected nilpotent Lie groups,
the group of Euclidean motions on the plane,
the Mautner group, or any closed  subgroup of one of these.
Our theorem gives the spectral invariance of $\GA$
if $A$ is the set of $C^{\infty}$-vectors
for the action of $G$ on $B$, or if $B= C_{0}(M)$, and
$A $ is  a set of $G$-differentiable Schwartz functions $\Cal S(M)$ on $M$.
This gives many examples of spectral invariant dense subalgebras
for the C*-algebras associated with dynamical systems.
We also obtain relevant
results about exact sequences, subalgebras, tensoring by
smooth compact operators, and strong spectral invariance in $L_{1}(G, B)$.
\vskip\baselineskip
\head Contents \endhead
$$\aligned &\text{   Introduction} \\
\S 1.&\text{   Strong Spectral Invariance}\\
\S 2.&\text{   Strong Spectral Invariance
of $\schwmh$ and $B^{\infty}$}\\
\S 3.&\text{   Subalgebras and Exact Sequences}\\
\S 4.&\text{   Tensoring by $M_{l}(\C)$ and Crossed Products
by Finite Groups}\\
\S 5.&\text{   Tensoring by Smooth Compacts}\\
\S 6.&\text{   Crossed Products by Type R Lie Groups and Spectral
Invariance
in $L_{1}(G, B)$}\\
\S 7.&\text{   Crossed Products by Polynomial Growth Groups}\\
&\text{   References.}
\endaligned
$$
\vskip\baselineskip
\heading  Introduction  \endheading
\par
The  theory of  differential
geometry on a noncommutative space
Connes \cite{6} requires
the use of \lq\lq differentiable structures\rq\rq\
for these noncommutative spaces,
or some sort of algebra of \lq\lq differentiable functions\rq\rq\
on the noncommutative space.  Such algebras of functions have usually been
provided by a dense subalgebra of smooth functions $A$
for which the $K$-theory $K_{*}(A)$ is the same as the
$K$-theory of the C*-algebra $K_{*}(B)$
(see for example Baum-Connes \cite{1}, Blackadar-Cuntz \cite{3},
Bost \cite{4}, Ji \cite{9}, the recent works of  G. Elliott, T. Natsume,
R. Nest,   P. Jollissaint, V. Nistor and many others).
In this paper, we use the algebras constructed in Schweitzer \cite{22} to
provide such dense subalgebras
for a large class of examples.
Let  $G$ be any compactly generated polynomial growth Type R Lie group,
not necessarily connected.
Here Type R means that all the eigenvalues of $Ad$ lie on the unit circle,
and polynomial growth means that the Haar measure of $U^{n}$ is bounded
by a polynomial in $n$, where $U$ is a generating neighborhood.
For example $G$ can be a
a finitely generated polynomial growth discrete group,
a compact or a connected nilpotent
Lie group, or the group of Euclidean motions on the plane,
the Mautner
group, or any closed  subgroup of one of these.
We provide
smooth subalgebras $G \rtimes \Cal S(M)$ of the C*-crossed product
$G\rtimes C_{0}(M)$, where $M$ is any $G$-space (see Examples 2.6-7,
6.26-7, 7.20).
We show that our subalgebras $G\rtimes \schwmns$ are all spectral invariant
in the C*-crossed product $G\rtimes C_{0}(M)$ (see Corollary
7.16), which implies that they have the
same $K$-theory as the C*-crossed product by
\cite{5}, VI.3
and \cite{21}, Lemma 1.2, Corollary 2.3.
For an example, if $M= H/K$ is a quotient of a compactly generated
polynomial growth Type R Lie group $H$ by
a closed cocompact subgroup $K$, and $G$ is a closed subgroup of
$H$, then we have the spectral invariance of $G\rtimes C^{\infty}(M)$.
If $G= \Z$, then $G\rtimes C^{\infty}(M)$ is just the standard
Fr\'echet algebra of Schwartz functions from $\Z$ to $C^{\infty}(M)$,
which are special cases of the algebras studied in Nest \cite{14}.
\par
We also show that if $B$ is any C*-algebra
on which $G$ acts, then the smooth crossed product
$G\rtimes B^{\infty}$ is spectral invariant in the C*-crossed
product $G\rtimes B$, where $B^{\infty}$ is the $C^{\infty}$-vectors for the
action of $G$ on $B$.  This generalizes the result Bost \cite{4},
Theorem 2.3.3
for elementary Abelian groups.
\par
To prove the spectral invariance of these dense subalgebras, we introduce
the notion of strong spectral invariance,
which implies
spectral invariance (see Definition 1.2 below -
this notion is similar to the condition (1.4) of Blackadar and
Cuntz \cite{3}, 3.1(b)).
 We show that
$\schwmns$ is
always strongly spectral invariant in $C_{0}(M)$,
and similarly for $B^{\infty}$ when $B$ is any Banach algebra with
an action of a Lie group.
If $A$ is
strongly spectral invariant in $B$,
we show the smooth
crossed product $\GA$, which consists of functions that
vanish rapidly  with respect to a generalized
\lq\lq word length function\rq\rq, is
strongly spectral invariant in $L_{1}(G, B)$. (For
$\GA$ to actually be a Fr\'echet algebra, we require that
 $G$  be compactly generated and Type R.
We also make the third assumption that
$G/{\text{Ker}}(Ad)$ has a cocompact solvable subgroup
 (this happens, for example, if $G$ is solvable or discrete).
We
conjecture that this third assumption is unnecessary - see
\cite{22}, \S 1.4,
Question 1.4.7.)
We then imitate a result of Pytlik \cite{19}
for the group algebra case to show
 that if $G$ has polynomial growth,
then a certain algebra of weighted $L_{1}$ functions $L_{1}^{\tau}(G, B)$
is in fact spectral invariant in the C*-crossed product $\GB$.  This
implies that $\GA$ is spectral invariant in $\GB$ - see
Corollaries 7.14 and
7.16. (For such polynomial growth groups we do not need the third
assumption above
about cocompact solvable subgroups, since such subgroups are
always present.
Hence our final result holds for an arbitrary compactly generated
polynomial growth Type R Lie group.)
\par
As a partial converse, for an arbitrary
Lie group $G$, if $G$ is {\it not} Type R,
then {\it none}
of the smooth Fr\'echet *-algebras
$\Cal S(G)$ we defined in \cite{22} are spectral invariant
in $L_{1}(G)$,
or in either of the C*-algebras
$C_{r}^{*}(G)$ or
$C^{*}(G)$ (see Theorem 6.29).
\par
We  analyze how the properties of spectral invariance and strong spectral
invariance behave with exact sequences, tensoring with $n\times n$ matrices
over $\C$ (see also \cite{21}), and tensoring by a
smooth version $\Kinf$ of the compact operators
$\K$.
For example, we show that the completed projective
tensor product $\Kinf \tensp A$ is spectral invariant in
the C*-algebra tensor product
$\K {\overline{\otimes}} B$ if $A$ is strongly spectral invariant in $B$.
Related results are obtained in Phillips \cite{18}, \S 4 for
the case $A=B$.
\par
The property of spectral invariance is important
for the study of how the representation theory of the dense subalgebra
relates to the representation theory of the C*-algebra.  For example,
$A$ is spectral invariant in $B$ iff every simple $A$-module
is contained in a *-representation of $B$ on a Hilbert space
(Schweitzer \cite{21}, Corollary 1.5, Lemma 1.2).
Our results on the spectral invariance of crossed products
by polynomial growth Lie groups thus generalize
the result of J. Ludwig \cite{12}  on the algebraically irreducible
representations of the Schwartz algebra of a nilpotent Lie group.
\par
Throughout this paper, the notations $\N$, $\Z$, $\R$, $\T$
shall be used for the natural numbers (with zero), integers, reals,
and the circle group respectively.  All of our algebras
will be over $\C$.  The term norm may be used
interchangably with the term seminorm.
If the positive definiteness
of a norm is  important, we shall state it explicitly.
The term differentiable will always mean infinitely
differentiable.
All groups will be assumed locally compact and Hausdorff.
\par
I would like to thank Chris Phillips for helpful comments and suggestions.
\vskip\baselineskip
\heading  \S 1 Strong Spectral Invariance  \endheading
\par
In the section, we define what it means for a dense Fr\'echet
subalgebra $A$ of a Banach algebra $B$ to be strongly spectral
invariant in $B$.
We show that strong spectral invariance implies spectral invariance,
and also exhibit an example of a spectral invariant dense subalgebra
which is not strongly spectral invariant.
\subheading{Definition 1.1}
A topological algebra is a topological vector
space over $\C$ with an algebra structure for which the multiplication
is separately continuous.
Let $A$ be a dense subalgebra of a topological algebra $B$.
If $A$ has no unit, let $\tilde A$ be $A$ with unit adjoined,
and let $\tilde B$ be $B$ with the same unit adjoined
(even if $B$ is unital already, adjoin a new one).
If $A$ has a unit, then $B$ has the same unit, and we let
${\tilde A}=A$, ${\tilde B}= B$.
We say that $B$ is a {\it $Q$-algebra} if $\tilde B$ has an open
group of invertible elements.
We shall usually be assuming that $B$ is a $Q$-algebra.
\par
By a {\it Fr\'echet algebra}, we mean a (locally convex) Fr\'echet space
with an algebra structure for which multiplication is jointly continuous
Waelbroeck \cite{24}, Chap VII.  A Fr\'echet algebra $A$ is {\it $m$-convex}
if there exists a family of submultiplicative seminorms on $A$
which give the topology of $A$.
We say that a Fr\'echet algebra
$A$ is a {\it dense Fr\'echet subalgebra} of
$B$ if $A$ is a dense subalgebra of $B$ and
the inclusion map $\iota \colon A\hookrightarrow B$
is a continuous injective algebra homomorphism.
Note that if $A$ is a dense (Fr\'echet) subalgebra
of  $B$, then $\tilde A$ is a dense (Fr\'echet)
subalgebra of   $\tilde B$.
\par
If $A$ is any dense subalgebra of
$B$, we say that ${A}$ is  {\it spectral invariant}
in $B$ if the invertible elements of $\tilde A$
are precisely those elements of $\tilde A$ which are invertible
in $\tilde B$ \cite{21}.
\subheading{Definition 1.2} Let $A$ be a dense Fr\'echet subalgebra
of a Banach algebra $B$.  Let $\bigl\{\pa \quad \pa_{m}\bigr\}$ be a
family of seminorms giving the topology of $A$, and arrange that
$\pa \quad \pa_{0}$
is a norm giving the topology of $B$. (From now on we shall
always assume that $\pa \quad \pa_{0}$ is a norm
giving the
topology on $B$ (though not always the
same one), whatever the choice of seminorms $\bigl\{
\pa \quad \pa_{m}\bigr\}$
topologizing $A$.)
We say that $A$ is {\it strongly spectral invariant } in $B$ if
$$\split (\exists C>0)(\forall m) & (\exists D_{m}>0)(\exists p_{m}\geq m)
(\forall n)(\forall a_{1},\dots a_{n} \in A)\\
& \biggl\{ \pa a_{1} \dots a_{n} \pa_{m}
\leq D_{m}C^{n}
\sum_{k_{1}+ \dots k_{n} \leq p_{m}}
 \pa a_{1} \pa_{k_{1}}\dots  \pa a_{n} \pa_{k_{n}}
\biggr\}.\endsplit \tag 1.3 $$
Notice that in the summand of (1.3), at most $p_{m}$ of the natural numbers
$k_{j}$
are nonzero, regardless of $n$.
This condition  appears similar to the condition
\cite{22}, (3.1.5) for  $m$-convexity, and also to \cite{22}, (3.1.19).
We require that $A$ be dense
in $B$ in order to show that strong spectral invariance implies
spectral invariance.   We choose to have
\lq\lq$\sum_{k_{1} + \dots k_{n} \leq p}$\rq\rq\
on the right hand side of (1.3) instead of
\lq\lq$\max_{k_{1} + \dots k_{n} \leq p}$\rq\rq\
since  sums commute with integration (see Theorems  5.4 and 6.7 below).
However both ways are equivalent up to a constant.
\par
We say that $A$ satisfies the {\it Blackadar-Cuntz   condition }
in $B$
if
there exists a family of seminorms
$\bigl\{ \pa \quad \pa_{m} \bigr\}$ for $A$
such that
$$ (\exists C>0)(\forall m)(\forall a,b\in A)\quad\biggl\{ \pa ab \pa_{m}
\leq C\sum_{i+j = m}
\pa a \pa_{i} \pa b \pa_{j} \biggr\} \tag 1.4 $$
(see Blackadar-Cuntz \cite{3}). Note that $m$ appears on the right hand
side of the inequality (1.4), whereas in (1.3) we replaced $m$
by the possibly larger natural number $p_{m}$.
\subheading{Example 1.5}  Let $B$ be the commutative
C*-algebra $c_{0}(\Z)$ of complex valued sequences on $\Z$
which vanish at infinity, with pointwise multiplication.  Let $A$ be
the dense Fr\'echet subalgebra $\Cal S (\Z)$ of sequences
which satisfy
$$ \pa f\pa_{m} = \sup_{n\in \Z} (1+|n|)^{m}\,|\,f(n)\,| <\infty,
\qquad m \in \N.\tag 1.6$$
Then  we have $\pa f_{1}\dots f_{n} \pa_{m} \leq
\pa f_{1}\pa_{m} \pa f_{2}\pa_{0}\dots \pa
f_{n}\pa_{0}$ for $f_{1}, \dots f_{n} \in A$.
So $A$
is strongly spectral invariant in $B$ (with $C=1$, $p_{m}=m$, $D_{m}=1$),
and also satisfies  the Blackadar-Cuntz condition in $B$, taking
$n=2$.
\proclaim{Proposition 1.7} Let $A$ be a dense Fr\'echet subalgebra
of a Banach algebra $B$.
Then $A$
is strongly spectral invariant in $B$ iff (1.3) holds for every
family of  seminorms $\bigl\{ \pa \quad \pa_{m} \bigr\}$ on $A$, or
iff (1.3) holds for any one family.
The constant $C$ depends only on the choice of the zeroth
norm $\pa \quad \pa_{0}$ on $B$.
If $A$ is strongly spectral invariant
in some Banach algebra $B$, then $A$ is $m$-convex.
\par
If $A$
satisfies the Blackadar-Cuntz condition in $B$, then
$A$ is strongly spectral invariant in $B$.  Moreover the constants
$D_{m}$ in (1.3) need not depend on $m$, and we may take $p_{m} = m$.
\endproclaim
\demo{Proof}
We first show that if (1.3) holds for one family of seminorms $\bigl\{
\pa \quad \pa_{m} \bigr\}$
on $A$, then it holds for any equivalent increasing family
$\bigl\{ \pa \quad\pa^{\prime}_{m} \bigr\}$ of seminorms on $A$.
First  we have
$$\aligned
\pa a_{1} \dots a_{n} \pa_{m}^{\prime} & \leq K\sum_{i\leq r }\biggl\{\pa
a_{1} \dots a_{n} \pa_{i} \biggr\}\qquad \qquad
{\text{equiv. of seminorms}}\\
&\leq K\sum_{i \leq r}\biggl\{
D_{i}C^{n} \sum_{k_{1} + \dots k_{n} \leq p_{i}}
\pa a_{1} \pa_{k_{1}} \dots \pa a_{n}
\pa_{k_{n}}\biggr\}, \endaligned \tag 1.8$$
where $r\in \N$ and $K>0$ depend only on $m$.
Let $t = \max_{i \leq r} p_{i}$.
Let ${ K_{1}}>0$ and $ s\geq t$
be such that for any $0\leq j \leq t$, we have
$\pa a \pa_{j} \leq { K_{1}} \pa a
\pa_{s}^{\prime}$.  (Here we have used that
the family $\bigl\{ \pa \quad \pa'_{m} \bigr\}$ is increasing.)
Define
$$ {\tilde {k_{j}}} = \cases 0 & k_{j}=0 \\ s & k_{j} \not= 0 \endcases$$
Then
the right hand side of (1.8) is bounded by
$$ \split K
\sum_{i \leq r} \biggl\{
D_{i}C^{n} \sum_{k_{1} + \dots k_{n} \leq t}
{ K_{1}}^{n}\pa a_{1} \pa^{\prime}_{\tilde{k_{1}}}
\dots \pa a_{n}& \pa_{
\tilde {k_{n}}}^{\prime}\biggr\}\\
\leq
&{ { {K_{2}}}}(CK_{1})^{n} \sum_{{ {k_{1}}} + \dots { {k_{n}}}
\leq st}
\pa a_{1} \pa_{{ {k_{1}}}}^{\prime} \dots \pa a_{n}
\pa^{\prime}_{ {k_{n}}}, \endsplit \tag 1.9 $$
where $K_{2} = K  \sum_{i \leq r} D_{i}$.
This shows that the increasing
family $\bigl\{ \pa \quad \pa_{m}^{\prime} \bigr\}$
satisfies (1.3).
A slightly longer argument shows that (1.3) holds for {\it any}
family of seminorms for $A$, but we omit it for brevity.
\par
To see that $C$ depends only on the
choice of the zeroth seminorm,
assume that in our above calculations that
$\pa \quad \pa_{0} = \pa \quad \pa^{\prime}_{0}$.  Then in (1.9) we could
replace $K_{1}^{n}$ with $K_{1}^{t}$ on the left hand side.
We would then have $C^{n}$ on the right hand side instead of
$(CK_{1})^{n}$, and set $K_{2} = KK_{1}^{t} \sum_{i\leq r} D_{i}$.
\par
For the $m$-convexity, it suffices to show that condition (1.3)
implies \cite{22}, (3.1.5) for every increasing family $\bigl\{
\pa \quad \pa_{m} \bigr\}$ of seminorms for $A$.  But by (1.3),
$$
\split
 \pa a_{1} \dots a_{n} \pa_{m}
& \leq D_{m}C^{n}
\sum_{k_{1}+ \dots k_{n} \leq p}
 \pa a_{1} \pa_{p}\dots  \pa a_{n} \pa_{p}
\\
& \leq {C_{1}}^{n}
 \pa a_{1} \pa_{p}\dots  \pa a_{n} \pa_{p}, \endsplit
 \tag 1.10 $$
where $C_{1}>0$
is a sufficiently large constant.
This is precisely \cite{22}, (3.1.5).
\par
Assume that $A$ satisfies the Blackadar-Cuntz condition in $B$.
Let $\bigr\{ \pa \quad \pa_{m} \bigl\}$ be a family of seminorms
which satisfy  (1.4).  Then we have
$$ \aligned
\pa a_{1} \dots a_{n} \pa_{m} & \leq
C \sum_{k_{1} + j = m } \pa a_{1} \pa_{k_{1}} \pa a_{2}
\dots a_{n}\pa_{j} \\
& \leq C^{2} \sum_{k_{1} + k_{2} + j = m}
\pa a_{1} \pa_{k_{1}} \pa a_{2}\pa_{k_{2}} \pa a_{3}
\dots a_{n}\pa_{j} \\
& \dots \qquad \dots \qquad \dots \qquad \dots \\
& \leq C^{n-1} \sum_{k_{1} + \dots k_{n} = m}
\pa a_{1} \pa_{k_{1}}  \dots\pa a_{n}\pa_{k_{n}}. \endaligned\tag 1.11$$
Taking $D_{m} = 1/C$, we have the last statement of the theorem.
\qed \enddemo
Motivated by \cite{22}, Theorem 3.1.4 and \cite{22},
Theorem 3.1.18, we ask the
following question.
\subheading{Question 1.12} If $A$ is strongly spectral invariant
in $B$, then
does $A$ satisfy the Blackadar-Cuntz condition (or some appropriate
modification of it) in $B$  ?
\vskip\baselineskip
\subheading{Example 1.13} We give an example of a spectral invariant
dense subalgebra which is not strongly spectral invariant.  Let $B $
be the C*-algebra $C^{*}(\Z)$
of the integers $\Z$, with convolution multiplication,
and let $A$ be the dense Banach
subalgebra $l_{1}(\Z)$ of absolutely summable
sequences.  Let $\pa \quad \pa_{0}$ be the C* norm,
and let $\pa \quad \pa_{1}$ be the $l_{1}$ norm.
We show that (1.3) cannot be satisfied for the
family of norms $\bigl\{ \pa \quad \pa_{0}, \pa \quad \pa_{1}
\bigr\}$ for $A$.
By the estimate \cite{10}, Chap VI, \S 6, p.82(in the
proof of Katznelson's theorem), we have
$$ \sup_{\psi = \psi^{*}, \,\pa \psi \pa_{1} \leq r}
\pa exp({i\psi }) \pa_{1}
= e^{r}. \tag 1.14$$
If $\psi = \psi^{*}$, note that $\pa exp(i\psi) \pa_{0}$ is just the
sup norm of the Fourier transform $e^{i{\hat\psi}}$ in $C(\T)$, which
is $1$ since $\hat \psi $ is real valued.   Also $\pa exp({i\psi}) \pa_{1}
\leq e^{\pa \psi \pa_{1}}$.
We have
$$ \aligned \pa exp(i \psi) \pa_{1} = \pa exp(i \psi/n)^{n} \pa_{1}
& \leq DC^{n} \pa exp(i \psi/n ) \pa_{1}^{p}\qquad {\text{by (1.3)}}\\
&\leq
DC^{n} (e^{r/n})^{p}\qquad \qquad {\text{if $\pa \psi
\pa_{1} \leq r$,}}\endaligned
\tag 1.15 $$
for all $n \in \N$.   Hence by (1.14),
$$ e^{r} \leq DC^{n} e^{rp/n} \tag 1.16$$
for arbitrarily large $r \in \N$.  This is a contradiction if we fix $n$
larger than $p$.   Hence
 $A$ is not strongly spectral invariant
in $B$.
However, $A$ is spectral invariant in $B$ by Wiener's theorem.
\vskip\baselineskip
\proclaim{Theorem 1.17} If $A$ is strongly spectral invariant in
$B$,
then $A$ is spectral invariant in $B$.
\endproclaim
\demo{Proof}
\proclaim{Lemma 1.18} If $A$ is strongly spectral invariant
 in $B$,
then $\tilde A$ is strongly spectral invariant in $\tilde B$.
\endproclaim
\demo{Proof}
It suffices to consider the case when $\tilde A$ and $\tilde B$ are
the respective unitizations of $A$ and $B$.
If $\bigl\{ \pa \quad \pa_{m} \bigr\}$ is a family of seminorms
for $A$, we define seminorms $\pa \quad \pa^{\prime}_{m}$ for
$\tilde A$ by
$\pa a + \lambda 1 \pa^{\prime}_{0}= \pa a \pa_{0} + | \lambda |$
and $\pa a + \lambda 1 \pa_{m}^{\prime} = \pa a \pa_{m}$ for $m >0$.
Let ${\tilde a}_{1}, \dots {\tilde a}_{n} \in {\tilde A}$,
where ${\tilde a}_{i} = a_{i} + \lambda_{i}1$, with
$a_{i} \in A$, $\lambda_{i} \in \C$.
We estimate $\pa {\tilde a}_{1}, \dots {\tilde a}_{n} \pa_{m}^{\prime}$.
For $m= 0$, the norm is submultiplicative up to a constant, so we
assume $m>0$.
Then
$$\aligned \pa {\tilde a}_{1}, \dots {\tilde a}_{n} \pa_{m}^{\prime}
&= \pa
({a}_{1} + \lambda_{1}1) \dots ({a}_{n} + \lambda_{n}1) \pa_{m}
\\
&\leq\sum_{1\leq i_{1} < \dots < i_{r}\leq n}
\pa a_{i_{1}}\dots a_{i_{r}}
\lambda_{j_{1}}\dots\lambda_{j_{n-r}} \pa_{m}, \endaligned  \tag 1.19 $$
where the sum is over all
$r$-tuples $1 \leq i_{1} < \dots i_{r} \leq n$ for all $1<r\leq n$,
and $\bigl\{ j_{1}, \dots j_{n-r}\bigr\} = \bigl\{ 1, \dots n\bigr\}
- \bigl\{ i_{1}, \dots i_{r} \bigr\}$.
We estimate one of the summands in $(1.19)$. For
simplicity, we set $D_{m} = 1$
in (1.3).  We have
$$\aligned
& \pa a_{i_{1}}\dots a_{i_{r}}
\lambda_{j_{1}}\dots\lambda_{j_{n-r}} \pa_{m}=
|\lambda_{j_{1}}|
\dots |\lambda_{j_{n-r}}|
 \pa { a}_{i_{1}}, \dots { a}_{i_{r}} \pa_{m}
\\
& \leq \,|\lambda_{j_{1}}|
\dots |\lambda_{j_{n-r}}| C^{r}\biggl\{   \sum_{k_{1}+ \dots k_{r} \leq p}
\pa a_{i_{1}} \pa_{k_{1}} \dots \pa a_{i_{r}} \pa_{k_{r}}
\biggr\} \qquad \text{strong spec. inv. (1.3)} \\
& \leq\,\,\, \pa {\tilde a}_{j_{1}} \pa^{\prime}_{0}
\dots\pa {\tilde a}_{j_{n-r}} \pa^{\prime}_{0} C^{r}\biggl\{
 \sum_{k_{1}+ \dots k_{r} \leq p}
 \pa {\tilde a}_{i_{1}} \pa^{\prime}_{k_{1}}
\dots \pa {\tilde a}_{i_{r}} \pa^{\prime}_{k_{r}}
\biggr\} \quad \text{def. of $\pa \quad \pa^{\prime}_{m}$.} \\
& \leq \, C^{r} \biggl\{ \sum_{k_{1}+ \dots k_{n} \leq p}
 \pa {\tilde a}_{1} \pa^{\prime}_{k_{1}}
\dots \pa {\tilde a}_{n} \pa^{\prime}_{k_{n}}
\biggr\}   \\
& \leq \, {\max(1, C)^{n}} \biggl\{\sum_{k_{1}+ \dots k_{n} \leq p}
  \pa {\tilde a}_{1} \pa^{\prime}_{k_{1}}
\dots \pa {\tilde a}_{n} \pa^{\prime}_{k_{n}}
\biggr\}  \text{} \endaligned \tag 1.20$$
Combining (1.19) and
(1.20), we have
$$ \pa {\tilde a}_{1} \dots {\tilde a}_{n} \pa_{m} \leq
2^{n}\max(1, C)^{n} \biggl\{ \sum_{k_{1}+ \dots k_{n} \leq p}
 \pa {\tilde a}_{{1}} \pa_{k_{1}}\dots \pa {\tilde a}_{n} \pa_{k_{n}}
\biggr\}. $$
This proves Lemma 1.18.
\qed \enddemo
\par
Now we are ready to prove Theorem 1.17.  By Lemma 1.18,
it suffices to consider the case when
$A$ and $B$ are both unital with the same unit.   By Proposition 1.7,
we may assume that the seminorms on $A$ are increasing
and submultiplicative. By strong spectral
invariance, we have
$$ \aligned \pa a^{n} \pa_{m} & \leq  D_{m}C^{n} \sum_{k_{1} + \dots k_{n}
\leq p}  \pa a \pa_{k_{1}} \dots \pa a \pa_{k_{n}} \\
& \leq D_{m}C^{n} \sum_{k_{1} + \dots k_{n}
\leq p}  \pa a \pa_{0}^{n-p}  \pa a \pa_{p}^{p} \quad
\text{since norms increase}\\
& \leq D_{m}^{\prime}(pC)^{n}  \pa a \pa_{0}^{n-p}
\pa a \pa_{p}^{p} \qquad
\endaligned \tag 1.21$$
where $p$
depends only on $m$, and $D_{m}^{\prime}$ is sufficiently large that
$D_{m} n^{p} \leq D_{m}^{\prime} p^{n}$ (note that $p$ is
fixed as $n$ runs).
If $\pa a\pa_{0}
< 1/2pC$, then by (1.21),
$$ \aligned \pa a^{n} \pa_{m} & \leq D_{m}^{\prime}(pC)^{n}
\bigl(1/2pC\bigr)^{n-p}
\pa a \pa_{p}^{p} \\
& = 1/2^{n} \biggl( (pC)^{p}D_{m}^{\prime} \pa a \pa^{p}_{p}2^{p} \biggr).
\endaligned  \tag 1.22 $$
So the series $1 + a + a^{2} + \dots $ converges absolutely in
the norm $\pa \quad \pa_{m}$ if $\pa a \pa_{0}$ is sufficiently small.
\par
Let $A_{m}$ be  the completion of $A$ in
$\pa \quad \pa_{m}$.  By what we've just seen, for each $m$
there is a neighborhood $U_{m}$ of the identity in $B$ such that
if $a \in A\cap U_{m}$, then the series $1+ (1-a) + (1-a)^{2} + \dots$
converges in $A_{m}$.  This is just the series for $a^{-1}$,
so $a^{-1} \in A_{m}$.
\par
Let $a \in A$ and assume $a^{-1}\in B$.  We show $a^{-1} \in A$.
Since we have choosen submultiplicative seminorms $\pa \quad \pa_{m}$
for $A$, it suffices to show that $a$ is invertible in $A_{m}$ for
each $m$ by Micheal \cite{13}, Theorem 5.2 (c).
\par
The set $a^{-1}U_{m} \cap U_{m} a^{-1}$
is open and nonempty (contains $a^{-1}$) in $B$, and so contains an element
$a^{\prime}$ of $A$ since $A$ is dense in $B$.  Then
$ a a^{\prime}$ and $a^{\prime}a$ are both in $U_{m}$.
Since they also lie in $A$, the construction of the $U_{m}$'s tells
us that $(aa^{\prime})^{-1}$ and $(a^{\prime}a)^{-1}$
both lie in $A_{m}$.  It follows that $a^{-1} \in A_{m}$.
Thus $a^{-1} \in A_{m} $ for all $m$ and  we are done.
The last part of this  argument is similar to Bost \cite{4}, Lemme A.2.3.
\qed
\enddemo
\heading  \S 2 Strong Spectral Invariance
 of $\schwmh$ and $B^{\infty}$\endheading
\par
We verify the strong spectral invariance of a space of Schwartz functions
$\schwmh$
on a locally compact $H$-space $M$ as a subalgebra of
the commutative C*-algebra
$C_{0}(M)$ of continuous functions vanishing at infinity on $M$.  Also,
we show that the set of $C^{\infty}$-vectors
$B^{\infty}$ is always strongly spectral invariant in $B$, for an arbitrary
Banach algebra $B$.
\par
We recall the definition of  $\schwmh$ from \cite{22}, \S 5.
 Let $H$ be a Lie group, possibly disconnected,
and let $M$ be  a locally compact space
on which $H$ acts.  We say that a Borel measurable function
$\s \colon M \rightarrow [0, \infty)$ is a {\it scale} if it is
bounded on compact subsets of $M$.
We say that a scale $\s$ {\it dominates} another scale $\gamma$
if there exists $C>0$ and $d\in \N$ such that
$ \gamma(m) \leq C(1+\s(m))^{d}$ for  $ m \in M$.
We say that $\s$ and $\gamma$ are {\it equivalent} (denoted
by $\s \thicksim \gamma$) if they dominate each
other.
If $h \in H$, define
$\s_{h}(m) = \s(h^{-1}m)$.   We say that $\s$ is {\it uniformly
$H$-translationally equivalent} if
for every compact subset $K$ of $H$ there exists $C_{K}>0$ and $d \in\N$
such that
$$ \s_{h} (m) \leq C_{K}(1+\s(m))^{d},\qquad m
\in M,\, h \in K. \tag 2.1 $$
If $\s$ is a uniformly $H$-translationally equivalent
scale on $M$,
we may define the $H$-differentiable $\s$-rapidly vanishing
functions $\schwmh$ by
$$ \schwmh = \{ f \in C_{0}(M), \,
\text{$f$ $H$-differentiable} \, | \,  X^{\gamma}f \in C_{0}(M)
\text{ and }
\pa \s^{d}X^{\gamma}f \pa_{\infty}
<\infty\,
\},  $$
where $X^{\gamma}$ ranges over all
 differential operators from the Lie algebra of
$H$, and $d$ ranges over all natural numbers.
We topologize $\schwmh$ by the  seminorms
$$ \pa f \pa_{d} = \sum_{|\gamma|= d}
\pa (1+\s)^{d}X^{\gamma}f \pa_{\infty}, $$
where we make the convention $|\gamma| = \sum |\gamma_{i}|$.
Then $\schwmh$ is an $m$-convex Fr\'echet *-algebra
under pointwise multiplication,
with differentiable action of $H$, and
is a dense Fr\'echet subalgebra of $C_{0}(M)$ \cite{22}, \S 5.
\par
In the following theorem, $B$ will be any Banach algebra, with a
strongly continuous action of a Lie group $G$ by isometric automorphisms.
(Throughout this paper, we
shall  assume that group actions on $B$ are   by isometric
automorphisms, which means that $\pa \alpha_{g}(b) \pa = \pa b \pa $
for all $b \in B$ and $g \in G$.)
We then may form
a dense Fr\'echet subalgebra $B^{\infty}$ of $C^{\infty}$-vectors
for the action of  $G$ on  $B$ \cite{22}, Theorem A.2.
\proclaim{Theorem 2.2} The Fr\'echet algebra $\schwmh$ is
a strongly spectral invariant subalgebra of $C_{0}(M)$,
and moreover satisfies the Blackadar-Cuntz condition in $C_{0}(M)$.
The same is true for the Fr\'echet algebra $B^{\infty}$
of $C^{\infty}$-vectors in $B$.
\endproclaim
\demo{Proof}
By the product rule,
$$\aligned \sum_{|\gamma|=d}
\pa (1+\s)^{d}X^{\gamma}(f_{1}  f_{2})\pa_{\infty}
&\leq C \sum_{|\gamma|=d,\,\beta_{1} + \beta_{2} = \gamma}
\pa (1+\s)^{d}(X^{\beta_{1}}f_{1})  (X^{\beta_{2}}f_{2}) \pa_{\infty}
\\
&\leq C\sum_{|\beta_{1}| +  |\beta_{2}| = d }
\pa (1+\s)^{d_{1}}(X^{\beta_{1}}f_{1}) \pa_{\infty}
\pa (1+\s)^{d_{2}}(X^{\beta_{2}}f_{2}) \pa_{\infty},
\endaligned
\tag 2.3 $$
where $d_{1}=|\beta_{1}|$ and $d_{2}=|\beta_{2}|$.
Since the right hand side of (2.3) is bounded by
$$ C \sum_{d_{1} + d_{2} = d}
\pa f_{1} \pa_{ d_{1}}
\pa f_{2} \pa_{d_{2}}, \tag 2.4 $$
we see that $\schwmh$ satisfies the
Blackadar-Cuntz condition in $C_{0}(M)$.
\par
We show that $B^{\infty}$ satisfies the
Blackadar-Cuntz condition in $B$.
Define seminorms by
$\pa b \pa_{m} = \sum_{|\gamma|= m}\pa X^{\gamma} b \pa$.
Then
$$ \sum_{|\gamma|=m}
\pa X^{\gamma}(b_{1}b_{2})\pa \leq C\sum_{|\gamma|=m,\,
\beta_{1} + \beta_{2}= \gamma}
\pa X^{\beta_{1}} b_{1}\pa \, \pa X^{\beta_{2}} b_{2} \pa =
C\sum_{i_{1}+i_{2} = m} \pa b_{1} \pa_{i_{1}} \pa b_{2} \pa_{i_{2}}.$$
\qed
\enddemo
\subheading{Definition 2.5}
We say that a locally compact group $H$ is {\it compactly generated}
if there
exists a relatively compact
neighborhood $U$ of the identity of $H$ such that
$U^{-1}= U$ and $H = \cup_{n=0}^{\infty}
U^{n}$.  We say that $H$ has {\it polynomial growth} if the Haar
measure of $U^{n}$ is bounded by a polynomial in $n$.  We say that
a scale $\tau$ on $H$ is a {\it gauge} if $\tau(e) = 0$, $\tau (g^{-1})=
\tau (g)$, and $\tau(gh) \leq \tau(g) + \tau (h)$.  We define the
{\it word gauge} $\tau_{U}$ on $H$ by
$$ \tau_{U}(g)= \min \{ \,n\,| \, g \in U^{n}\,\}.$$
Then $\tau_{U}$ is independent up to equivalence
of the choice of $U$ \cite{22}, Theorem 1.1.21.
\subheading{Example 2.6}  Assume that $H$ is a compactly
generated polynomial
growth Lie group, and let $K$ be a closed subgroup.  Let $\s(h) =
\inf_{k \in K}
\tau_{U}(hk)$,
and $M= H/K$. Then $\Cal  S_{H}^{\s}(M)$ is strongly
spectral invariant in
$C_{0}(M)$.  Moreover, $\Cal S_{H}^{\s}(M)$ is also a nuclear
Fr\'echet algebra \cite{22}, Proposition 1.5.1, Theorem 6.8.
\subheading{Example 2.7}
Let $G$ be any other closed subgroup of $H$.  Then $\Cal S_{G}^{\s}(M)$
is also strongly spectral invariant in $C_{0}(M)$.    However
$\Cal S_{G}^{\s}(M)$ may not be nuclear. For example, let $G=K=\{ e\}$,
$H= \T$.  Then $\Cal S_{G}^{\s} (M)$ is
the infinite dimensional Banach algebra $ C(\T)$
of continuous functions on the circle, and hence not nuclear.
\par
To give a familiar example, let $M=H= \R$, $G=K= \{ e \}$.  Then $\s $
is equivalent to the absolute value function on $\R$, and $\Cal S_{H}^{\s}
(M)$
is just the standard set of Schwartz functions $\schwr$ on $\R$.
The space $\Cal S_{G}^{\s}(M)$ is the set of continuous
functions on $\R$, which vanish rapidly with respect to $|r|$.
\heading \S 3 Subalgebras and Exact Sequences \endheading
We look at how the properties of spectral invariance and
strong spectral invariance behave in the context of exact sequences.
We let $I$ be a  two-sided ideal of $A$, and let $J$ be the closure
of $I$ in $B$.  Note that if we assume that $J \cap A = I$,
then $A/I$
is a subalgebra  $B/J$.
\subheading{Example 3.1} For example, for $A= C^{\infty}(\T)$, $B= C(\T)$,
we could take $I$ and  $J$ to be the set of functions in $A$ and $B$
respectively which vanish at some fixed point $p\in \T$.  The
property $J \cap A = I$ is then satisfied.  However, if in place of $I$
we took the ideal of functions which vanish along with all of their
derivatives at $p$, we would still have ${\overline I} = J$ but not
$J\cap A = I$.
\par
Results related to part (2) of the following theorem appear in Palmer
\cite{15},
Corollary 5.6,7 (see also the introduction of that paper).
(In that paper,
a \lq\lq spectral invariant subalgebra\rq\rq\
is called a \lq\lq spectral subalgebra\rq\rq.)
\proclaim{Theorem 3.2}
Let $I$ and $J$ be as above, and assume $J \cap A = I$.
Let
$A_{1}$ be any  subalgebra of $A$, and  let $B_{1}$ be the closure
of $A_{1}$ in $B$.  Assume that $B$ is a Banach algebra, so that
both $B_{1}$ and $B/J$ are Banach algebras and hence $Q$-algebras.
\roster
\item
Let $A$ be a dense Fr\'echet subalgebra of  $B$,
and assume that $I$ and $A_{1}$ are both closed in the topology of
$A$, with Fr\'echet topology inherited from $A$.
Let $A$ be strongly spectral invariant in $B$. Then $A_{1}$ is strongly
spectral invariant in $B_{1}$, and the ideal $I$
is strongly spectral  invariant
in $J$.
Similar statements hold for the Blackadar-Cuntz condition.
\item
Let $A$ be any dense subalgebra of $B$.
Assume that $A$ is spectral invariant in $B$.  If $A_{1}= A\cap B_{1}$,
then $A_{1}$ is
spectral invariant in $B_{1}$.  The ideal $I$ is  spectral  invariant
in $J$ and $A/I$ is  spectral invariant in $B/J$.
Conversely, if $I$ is spectral invariant in $J$ and $A/I$
is spectral invariant in $B/J$, then  $A$ is spectral invariant in $B$.
\endroster
\endproclaim
\demo{Proof}
Since  seminorms on $A_{1}$ are given by the restriction of any family
of seminorms on $A$ to $A_{1}$, and the norm on $B_{1}$ is the restriction
of the norm on $B$ to $B_{1}$, the strong spectral invariance of $A_{1}$
in $B_{1}$ is obtained simply by restricting the inequality (1.3)
in the definition  of strong spectral invariance
to elements of $A_{1}$.
As a special case, we also see that $I$ is strongly spectral invariant in
$J$. This proves (1).
\par
Assume that $A$ is spectral invariant in $B$.
Recall from \cite{21}, Theorem 1.4 that $A$ is spectral invariant
in $B$ iff every simple $A$-module is contained
as a dense $A$-submodule of a
continuous $B$-module.
Let $V$ be simple  $A/I$-module.
Then $V$ is a simple $A$-module and so extends to a continuous
$B$-module $W$ in which $V$ is dense (using
\cite{21}, Theorem 1.4 $(i)\Rightarrow
(iii)$).
To see that $W$ factors to
a $B/J$-module it suffices to show that $JW=\{0\}$.
We know $IV=\{0\}$.  So since $J=\overline{I}$ and
$V$ is dense in $W$, $JW=\{0\}$.
Thus $W$ is a $B/J$-module extending $V$, and we have
shown that $A/I$
is spectral invariant in $B/J$ by \cite{21}, Theorem 1.4 $(iii) \Rightarrow
(i)$. This did not require $B/J$ to be a $Q$-algebra, so it would
suffice for $B$ to be any $Q$-algebra.
\par
Next we show that
$A_{1}$ is spectral invariant in $B_{1}$
under the assumption that  $A_{1}= A\cap B_{1}$ (and hence,
as a special case,
that $I$ is spectral invariant in $J$).
By replacing $A$, $B$ with $\tilde A$, $\tilde B$,
we may assume that $A$ and $B$ are unital with the same unit.
If $A_{1}$ and $B_{1}$ are unital, let $a \in A_{1}$, $a^{-1} \in B_{1}$.
Then $aa^{-1}=a^{-1}a = 1_{B_{1}}$.  Let $q$ be the projection
$1_{B}-1_{B_{1}}$.  Then $(q + a)(q + a^{-1}) =(q+a^{-1})(q+a)= 1_{B}$,
so $q + a$ is invertible in $B$ and hence in $A$.
Hence $q+a^{-1} \in A$, and $a^{-1} \in A$.  By assumption,
$a^{-1} \in A_{1}$.
\par
Next assume $A_{1}$
is nonunital.  We may assume that ${\tilde A}_{1}$, ${\tilde B}_{1}$,
$A$, and $B$ all have the same unit.  If $a \in {\tilde A}_{1}$ is
invertible in ${\tilde B}_{1}$, then clearly $a$ is invertible in $B$.
Hence $a^{-1}$ lies in both $A$ and ${\tilde B}_{1}$, so by assumption
$a^{-1} \in {\tilde A}_{1}$ and $A_{1}$ is spectral invariant in $B_{1}$.
\par
Assume that
$I$ is spectral invariant in $J$ and $A/I$ is spectral invariant in $B/J$.
We show that $A$ is spectral invariant in $B$
using \cite{21}, Theorem 1.4.  Let $V$ be an irreducible
$A$-module.  One easily checks that either
$IV=\{0\}$ or $V$ is a simple $I$-module.
So we have two cases.
\par
Case 1.  Say $V$ is a simple $ I$-module.
Then since $I$ is spectral invariant in $J$, we can extend
$V$ to a $J$-module $W$.
Since $J$ is a two-sided ideal in $B$, we can
extend the action of
$J$ on $W$ in a unique way to one of
$B$ (see, for example, the argument in Fell \cite{7}, Proposition 1).
\par
Case 2.  Say $IV=\{0\}$.  Then $V$ is an irreducible
$A/I$-module.  By hypothesis, and since $B/J$ is a $Q$-algebra,
$V$ extends to a
$B/J$-module $W$.
Using the
canonical algebra homomorphism from $B$ to $B/J$, we make $W$ a
$B$-module.
This action of $B$ on $W$ clearly extends the action of
$A$ on $V$.
This proves (2).
\qed
\enddemo
\subheading{Remark 3.3} Theorem 3.2(2) and \cite{21}, Lemma 1.2 answer
Question 3.1.8 of Blackadar \cite{2} in the affirmative, in the case that
the local Banach algebra $A$ has a Fr\'echet topology stronger than
the norm topology, and the ideal $I$ is closed in $A$.
A general answer to
this question is given in Schmitt \cite{20}.
\heading \S 4
Tensoring by $M_{l}(\C)$ and Crossed Products by
finite groups \endheading
\proclaim{Theorem 4.1}
\roster
\item Let $A$ be a dense Fr\'echet subalgebra
of a Banach algebra $B$.
Then the matrix algebra $M_{l}(A)$ is strongly
spectral invariant in $M_{l}(B)$
iff $A$ is strongly spectral invariant in $B$.
\item (Schweitzer \cite{21}, Theorem 2.1) Let $A$ be a dense subalgebra
of a $Q$-algebra $B$.  Then the matrix algebra $M_{l}(A)$
is spectral invariant in $M_{l}(B)$ iff
$A$ is spectral invariant in $B$.
\endroster
\endproclaim
\demo{Proof}
For part (1), by Theorem 3.2 (1) it suffices to show that if $A$ is
strongly spectral invariant in $B$, then $M_{l}(A)$ is strongly
spectral invariant in $M_{l}(B)$.  We define seminorms $\pa \quad
\pa_{m}^{\prime}$
on $M_{l}(A)$ by
$$ \pa [a ]\pa_{m}^{\prime}
= \max_{1 \leq i, j \leq l } \pa [a]_{ij} \pa_{m}. $$
An inductive argument shows
that the $ij$th entry of a product $[a_{1}]\dots [a_{n}]$
is the sum of $l^{n-1}$ products of elements of the form
$[a_{1}]_{i_{1}j_{1}} \dots [a_{n}]_{i_{n} j_{n}}
$.
By the strong spectral invariance of $A$ in $B$, we have
$$ \aligned
\pa [a_{1}]_{i_{1}j_{1}}\dots [a_{n}]_{i_{n}j_{n}} \pa_{m}
&
\leq
D_{m}C^{n} \sum_{k_{1}+ \dots k_{n} \leq p}
\biggr\{ \pa [a_{1}]_{i_{1}j_{1}} \pa_{k_{1}}\dots
\pa [a_{n}]_{i_{n}j_{n}} \pa_{k_{n}}  \biggr\}
\\
& \leq
D_{m}C^{n} \sum_{k_{1}+ \dots k_{n} \leq p}
\biggr\{ \pa [a_{1}] \pa_{k_{1}}^{\prime}\dots
\pa [a_{n}] \pa^{\prime}_{k_{n}}  \biggr\}.
\endaligned
$$
It follows that
$$
\pa [a_{1}]\dots [a_{n}] \pa^{\prime}_{m}
\leq
D_{m}l^{n-1}C^{n} \sum_{k_{1}+ \dots k_{n} \leq p}
\biggr\{ \pa [a_{1}] \pa_{k_{1}}^{\prime}\dots
\pa [a_{n}] \pa^{\prime}_{k_{n}}  \biggr\}.
$$
\qed
\enddemo
\proclaim{Corollary 4.2} Let $G$ be a finite group acting on $A$ and $B$
by algebra automorphisms, which are continuous on $B$.
\roster
\item Let $A$ be a dense Fr\'echet subalgebra of a
Banach algebra $B$, and assume that $G$ acts continuously on $A$. Then the
crossed product $G\rtimes A$ is strongly spectral invariant
in $G\rtimes B$ iff $A$ is strongly spectral invariant in $B$.
\item Let $A$ be a dense subalgebra of a $Q$-algebra $B$.
Then the crossed product $G\rtimes A$ is spectral invariant
in $G\rtimes B$ iff $A$ is  spectral invariant in $B$.
\endroster
\endproclaim
\demo{Proof}
First assume that $G\rtimes A$ is (strongly) spectral invariant
in $G\rtimes B$.  Then by (Theorem 3.2 (1)) Theorem 3.2 (2), and since
$B\cap G \rtimes A =A$, we have that $A$ is (strongly)
spectral invariant in $B$.
\par
Next assume that $A$ is (strongly) spectral invariant in $B$.
The following argument is similar to one in the appendix of
Baum-Connes \cite{1}.
Let $l$ be the order of $G$.  We may identify $M_{l}(B)$ with
the set of functions from $G \times G $ to $B$, if we define the
multiplication
$$ S* T (g, h) = \sum_{k \in G} S(g, k)T(k, h) $$
on elements $S, T $ of $C(G\times G, B)$.
We make similar definitions for $M_{l}(A)$.  Then by Theorem 4.1,
we know that $C(G\times G, A)$ is (strongly) spectral invariant
in $C(G\times G, B)$.
\par
We embed $G\rtimes B$ in $C(G\times G, B)$ via
$$i(F)(g, h) = \alpha_{g}(F(g^{-1}h)), \qquad F \in G\rtimes B.  $$
This embedding is easily seen to be an algebra homomorphism, and
is a topological embedding since $G$ acts continuously.
Similarly $i$ embeds $G\rtimes A$ as a subalgebra of $C(G\times G, A)$,
topologically if $A$ is Fr\'echet.
\par
Let $G$ act on $C(G\times G, B)$ (resp. $C(G\times G, A)$) via
$$ \theta_{g}(S)(k, h) = \alpha_{g^{-1}}(S(gk, gh)).$$
Then $i(G\rtimes A)$ is the set of fixed points for the
action of $\theta$ on $C(G\times G, A)$, and similarly for
$i(G\rtimes B)$ and $C(G\times G, B)$.  Clearly
$i(G\rtimes A) = i(G\rtimes B) \cap C(G\times G, A)$.
So (Theorem 3.2 (1)) Theorem 3.2 (2)
tells us that $G\rtimes A$ is (strongly) spectral
invariant in $G\rtimes B$.
\qed
\enddemo
\subheading{Remark 4.3} There is a nice alternate proof of the spectral
invariance part of Corollary 4.2 using extensions of simple modules
and Fell \cite{7}, Proposition 5.
\heading \S 5 Tensoring  by Smooth Compacts
\endheading
\par
Throughout this section, we
let $B$ be a C*-algebra, and
$A\subseteq B$ be a dense Fr\'echet subalgebra.
Let $\Cal H = l_{2}(\Z )$ and let $\K$ be
the compact operators on $\H$.
We define two dense subalgebras of the C*-tensor
product $\K {\overline{\otimes}} B$.
\par
If $\Cal A$ is any Fr\'echet algebra, let $\Cal S(\Z^{2}, \Cal A)$
be the set of  functions $\varphi$ from $\Z^{2}$ to
$\Cal A$ which satisfy
$$
\pa \varphi
\pa_{q} = \sum_{r, s\in \Z} \bigl( 1 + |r| + |s|\bigr)^{q}
\pa \varphi(r, s) \pa_{q} <\infty\tag 5.1$$
for every $q \in \N$.
We define multiplication
by
$$ \varphi * \psi (r, t) = \sum_{s\in \Z^{2}} \varphi(r, s) \psi(s, t).
\tag 5.2$$
If $\Cal A = \C$, we denote the resulting nuclear
$m$-convex Fr\'echet algebra
by $\Kinf$, the
 {\it smooth compact operators}.
In general,
$\Cal S(\Z^{2}, \Cal A)$ is isomorphic to the projective
completion $\Kinf \tensp \Cal A$.
Let $\Cal L$ be a Hilbert space on which $B$ is faithfully
*-represented.  Then $\K {\overline {\otimes}} B$ is faithfully
represented on the Hilbert space tensor product
$\H \otimes \Cal L= l_{2}(\Z, \Cal L )$, and
$\Kinf \tensp A$ is the dense subalgebra of $\K {\overline{\otimes}} B$
which acts on $l_{2}(\Z, \Cal L)$ via
$$ \varphi\xi (r) = \sum_{t\in \Z} \varphi(r, t)\xi (t),
\qquad \varphi \in \Kinf
\tensp A, \quad \xi \in l_{2}(\Z, \Cal L). \tag 5.3$$
\par
We also note that the space
$l_{1}(\Z^{2}, B)\cong l_{1}(\Z^{2}) \tensp B$,
topologized by the norm  (5.1) with $q = 0$,
is a Banach algebra under the multiplication (5.2).  It also
acts on $l_{2}(\Z^{2}, \Cal L)$ via (5.3) and is a dense subalgebra
of $\K {\overline{\otimes}} B$, which contains $\Kinf \tensp A$
as a dense subalgebra.
\proclaim{Theorem 5.4} If $A$ is strongly spectral invariant in
$B$, then $\Kinf \tensp A$ is strongly spectral invariant in $l_{1}(\Z^{2})
\tensp B$.
\endproclaim
\demo{Proof}
Let $\varphi_{1},\dots \varphi_{n} \in \Kinf \tensp A$.
We estimate
$$ \aligned
&\pa \varphi_{1} * \dots \varphi_{n} \pa_{q}  \leq
\sum_{r_{1}, \dots r_{n+1} \in \Z}
(1+ |r_{1}| + |r_{n+1}|)^{q}
\pa \varphi_{1}(r_{1}, r_{2}) \dots \varphi_{n} (r_{n}, r_{n+1}) \pa_{q}
\\
&
\leq
D_{q}C^{n}\sum_{k_{1} + \dots k_{n} \leq p_{q}}\biggl\{
\sum_{r_{1}, \dots r_{n+1} \in \Z}
(1+ |r_{1}| + |r_{n+1}|)^{q}
\pa \varphi_{1}(r_{1}, r_{2})\pa_{k_{1}}
 \dots \varphi_{n} (r_{n}, r_{n+1}) \pa_{k_{n}}\biggr\}.
\endaligned \tag 5.5 $$
Because of the inequality (see argument after \cite{22}, (3.2.5))
$$(1 + |r_{1}| + |r_{n+1}|)^{q} \leq 2^{q}\bigl\{ (1 + |r_{1}|)^{q}
 +( 1 + |r_{n+1}|)^{q} \bigr\} ,$$
the right hand side of (5.5) is bounded by
$$ \split D_{q}2^{q} C^{n} \sum_{k_{1} +\dots k_{n} \leq p_{q}} &\biggl\{
\pa \varphi_{1}\pa_{k_{1}+ q} \dots \pa \varphi_{n} \pa_{k_{n}}+
\pa \varphi_{1}\pa_{k_{1}} \dots \pa \varphi_{n} \pa_{k_{n} + q} \biggr\}
\\
& \leq D_{q}2^{q}{ C}^{n} \sum_{k_{1}+ \dots k_{n} \leq p_{q} + q}
\biggl\{
\pa \varphi_{1} \pa_{k_{1}} \dots \pa \varphi_{n} \pa_{k_{n}} \biggr\}.
\endsplit$$
This gives the strong spectral invariance of $\Kinf \tensp A$
in $l_{1}(\Z^{2})\tensp B$.
\qed
\enddemo
\proclaim{Corollary 5.6}  If $A$ is strongly spectral invariant in
$B$, then $\Kinf \tensp A$ is spectral invariant in  $\K
{\overline{\otimes}} B$.\endproclaim
\demo{Proof}
By Theorem 5.4, it suffices to show that $C=l_{1}(\Z^{2}) \tensp B$ is
spectral invariant in $D=\K {\overline{\otimes }} B$.
First we note that $CDC \subseteq C$.  For let $\psi_{1},\psi_{2} \in
C$ and $\varphi \in D$.  Then, in the above representation
(5.3) on $l_{2}(\Z, \Cal L)$,  $\varphi$ may be thought of as
a $\Z \times \Z$ matrix with entries in $B$, where each entry
has norm bounded by $c=\pa \varphi \pa_{\K {\overline{\otimes}} B}$.
We have
$$ \split\pa\psi_{1} *\varphi *\psi_{2} \pa_{0}
& \leq \sum_{r, s, t,w} \pa \psi_{1}(r, s) \pa_{B} \pa
\varphi(s, t)\pa_{B}
\pa \psi_{2}(t, w) \pa_{B} \\
&
\leq
c \pa \psi_{1} \pa_{0} \pa \psi_{2} \pa_{0}, \endsplit $$
where $\pa \quad \pa_{0}$ is the norm on $C$ (see (5.1) with $q=0$).
So $\psi_{1} * \varphi * \psi_{2} \in C$ and $CDC \subseteq C$.
We see that $C$ is spectral invariant in $D$ by the following lemma.
\proclaim{Lemma 5.7 (Compare Bost \cite{4}, Proposition A.2.8)}
Let $A$ be any subalgebra of an algebra $B$.
If $ABA\subseteq A$, then $A$ is spectral invariant in $B$.
\endproclaim
\demo{Proof}
If $A$ is unital, then $A=B$ and we are done.  So assume that
$A$ is nonunital.
Let $a+\lambda 1 \in \Atil$ for $\lambda \not= 0$,
and assume that
$(a+ \lambda 1)^{-1} \in \tilde B$.
Then
there is a
$b\in B$ such that $b + 1/{\lambda} =(a + \lambda 1)^{-1}\in \Botil $.
Since $(a+\lambda 1 )(b+ {1/{\lambda}})
=(b+ {1/{\lambda}})(a+\lambda 1 ) = 1$, we have $b = {{-a/{\lambda^{2}}}
- {ab/\lambda  }} = {{-a/{\lambda^{2}}}
- {ba/{\lambda}  }}$.  Substituting the first into the
second  gives $b \in ABA \subseteq A$.
So $a + \lambda 1$ is invertible in $\tilde A$
and $A$ is spectral invariant in $B$.
\qed \enddemo
This proves Corollary 5.6.
\qed
\enddemo
\par
Lemma 5.7 can also be used to show directly that $\Cal C=
\Kinf \tensp B$ is spectral invariant
in $\Cal D=\K {\overline \otimes} B$, since $\Cal C \Cal D \Cal C
\subseteq
\Cal C$.
\subheading{Remark 5.8}
The proof of the spectral invariance part of Corollary
4.2 above is easily
generalized from crossed products of finite groups
to crossed products by compact Lie groups $G$,
if we assume that $A$ is strongly spectral invariant in $B$, and
use Corollary 5.6, with $G$ in place of $\Z$.  We omit the details,
since we will be obtaining this same result using other methods in
\S 7 (see Corollary 7.16).
\heading \S 6 Crossed Products by Type R Lie Groups and Spectral
Invariance
in $L^{1}(G, B)$ \endheading
\par
We show that strong spectral invariance is preserved by
taking crossed products, if we use $L_1(G, B)$ as the
Banach algebra crossed product, and use a subadditive scale
on $G$ to define the smooth crossed product $G\rtimes^{\tau} A$.
Throughout this section, $A$ will be a dense Fr\'echet subalgebra
of a Banach algebra $B$.
\subheading{Definition 6.1}
We recall some definitions from \cite{22}, \S 2.
Let $\tau$ be
a scale on a locally compact group $G$ (see \S 2 above).
Let $E$ be any Fr\'echet space.
We define the {\it $\tau$-rapidly vanishing $L_{1}$ functions
$L_{1}^{\tau}(G, E)$ from $G$ to $E$} to be the set of
measurable functions
$\varphi$ from $G$ to $E$ such that
$$\pa \varphi \pa_{d, m} =
\int_{G}\pa\tau^{d}\varphi(g) \pa_{m} dg <\infty, \tag 6.2$$
where $\pa \quad \pa_{m}$ ranges over a family of seminorms for $E$,
and $d$ ranges over the natural numbers.  We shall always use
$dg$ to designate a {\it left} Haar measure on $G$.
We topologize $L_{1}^{\tau}(G, E)$ by the seminorms (6.2).
\par
Topologize $E$ by an increasing family of seminorms $\bigl\{
\pa \quad \pa_{m} \bigr\}$.
We say that the action of $G$ on a $G$-module $E$ is $\tau$-tempered
if for every $m \in \N$ there exists a polynomial $\text{poly}$
and $k \in \N$ such that
$$ \pa \alpha_{g}(e)
\pa_{m} \leq \text{poly}(\tau(g)) \pa e \pa_{k}, \qquad
e \in E, \,g \in G.\tag 6.3
 $$
We say that $\tau$ is {\it sub-polynomial} if there exists a polynomial
$\text{poly}$ such that
$$ \tau(gh) \leq \text{poly}(\tau(g),\tau(h)),
\qquad g,\, h \in G. $$ From the estimates
of \cite{22}, Theorem 2.2.6,
we see that if $A$ is a Fr\'echet
algebra with $\tau$-tempered action of $G$, and $\tau$ is sub-polynomial,
{\it then $L_{1}^{\tau}(G, A)$ is a Fr\'echet algebra under convolution}.
\par
Next, we assume in addition that $G$ is a Lie group, possibly
disconnected.
 We define the {\it differentiable $\tau$-rapidly
vanishing functions $\Cal S_{1}^{\tau}(G, E)$ from $G$ to $E$}
to be the set
of differentiable functions $\varphi$ from $G$ to $E$ such that
$$\pa \varphi \pa_{d, \gamma, m} =
\int_{G}\pa\tau^{d}X^{\gamma}\varphi(g) \pa_{m} dg <\infty, \tag 6.4$$
where $X^{\gamma}$ is any differential operator from the Lie
algebra of $G$ acting by left translation, and $d$ and $\pa \quad \pa_{m}$
are as in (6.2).
We topologize $\Cal S_{1}^{\tau}(G, E)$ by the seminorms (6.4).
\par
To make our
definitions sufficiently general, we let $H$ be any Lie group containing
$G$ as a subgroup with differentiable inclusion map.
We say that $\tau$ {\it bounds $Ad$ on $H$} if there exists a polynomial
$\text{poly}$
such that
$$ \pa Ad_{g} \pa \leq \text{poly}(\tau(g)), \qquad g \in G, \tag 6.5  $$
where $\pa Ad_{g} \pa$ is the operator norm of $Ad_{g}$ as an operator
on the Lie algebra of $H$.
The {\it inverse scale} $\tau_{-}$ is defined by
$\tau_{-}(g) = \tau(g^{-1})$.
And finally, if $\tau$ is a sub-polynomial scale on $G$ such that
$\tau_{-}$ bounds $Ad$ on $H$, and $\s$ is  an $H$-translationally
equivalent scale on a locally compact $H$-space $M$ (see \S 2),
and we have
$$ \s(gm) \leq \text{poly}(\tau(g),\s(m)),
\qquad g\in G, m\in M,\tag 6.6$$
for some polynomial $\text{poly}$, then
we say that {\it $(M,\s, H)$ is a
scaled $(G, \tau)$-space}.
By \cite{22},
Theorem 2.2.6,  Theorem 5.17
it follows that if
$\tau$ is a sub-polynomial
scale on $G$ such that
the action of $G$ on the Fr\'echet algebra $A$ is $\tau$-tempered,
and  either
$\tau_{-}$ bounds $Ad$ on $G$ or $G$ acts differentiably on
$A$,
{\it then $\Cal S^{\tau}_{1}(G, A)$ is a Fr\'echet
algebra under convolution,
which we denote by $G\rtimes^{\tau} A$}.
{\it Moreover, if $(M, \s, H)$ is a scaled $(G, \tau)$-space,
then the action of $G$ on $\schwmh$ is
$\tau$-tempered.  In particular, $G \rtimes^{\tau}
\schwmh$ is a Fr\'echet
algebra}.
\par
We warn the reader that if $G$ acts differentiably on $A$, and the
action of $G$ on $A$ is $\tau$-tempered, it often
is a prerequisite that $\tau_{-}$
must
bound $Ad$ on $G$.
Hence we often gain nothing by having the either/or
hypothesis in the previous paragraph.
On the other hand, in the group algebra case, when $A= \C$,
there is no requirement on $\tau$ except that $\tau$ be sub-polynomial.
\par
In light of these comments,
we shall  often require the Lie group $G$ to have
a gauge that bounds $Ad$ (see  Definition 2.5 for
the definition of a gauge),
in order to apply the second and third paragraphs of Theorem 6.7 below.
We say that $G$ is {\it Type R} if for every $g \in G$,
$Ad_{g}$, as an operator on the Lie algebra of $G$, has eigenvalues
on the unit circle.
Assume that $G$ is compactly generated, and that the
group $G/Ker(Ad)$ has a
cocompact solvable subgroup (this holds, for example, when $G$ is
solvable or discrete).  Then {\it $G$ has a
gauge that bounds $Ad$ iff $G$ is Type R} \cite{22}, Theorem 1.4.3.
If $G$ is
Type R, then in fact the word gauge (see Definition 2.5) bounds $Ad$.
\par
Examples of Type R Lie groups are given by discrete groups,
and closed subgroups of connected polynomial
growth Lie groups.  See the introduction and
\cite{22}, \S 1.4 for more examples.
\proclaim{Theorem 6.7}  Let $G$ be a locally compact group and
let $\tau$ be a subadditive scale
on $G$. ($\tau$ is {\it subadditive } if $\tau(gh)
\leq \tau(g) + \tau(h)$.)
Assume that the action of  $G$ on $A$ is
$\tau$-tempered. Then the Fr\'echet algebra $L_{1}^{\tau}(G, A)$
is strongly spectral invariant
in $L_{1}(G, B)$ if $A$ is strongly spectral invariant in  $B$.
\par
If, in addition, $G$ is a Lie group, and either $\tau_{-}$ bounds $Ad$
on $G$ or $G$ acts differentiably
on $A$, then $G\rtimes^{\tau} A$ is
a Fr\'echet algebra which is strongly spectral invariant in
$L_{1}(G, B)$ if $A$ is strongly spectral invariant in $B$.
\par
In particular, if $(M, \s, H)$ is a scaled $(G, \tau)$-space, then
the smooth crossed product $G\rtimes^{\tau}\schwmh$ is strongly
spectral invariant
in  $L_{1}(G, C_{0}(M))$.  Also, $G\rtimes^{\tau}B^{\infty}$
is strongly spectral invariant in $L_{1}(G, B)$, where $B^{\infty}$
denotes the set of $C^{\infty}$-vectors for the action of $G$ on $B$.
($G$ must often be a Type R Lie group
- see remarks
preceding the theorem.)
\endproclaim
\vskip\baselineskip
\proclaim{Corollary 6.8}  If $G$ is any Lie group
and $\tau$ is any subadditive scale on $G$, then the group Schwartz
algebra $\Cal S_{1}^{\tau}(G)$ is strongly spectral invariant
in $L_{1}(G)$.  In particular, this holds if $\tau$ is the word
gauge  on a compactly generated Lie group $G$.
\endproclaim
\vskip\baselineskip
\demo{Proof} By Theorem 2.2,
it suffices to prove the first two paragraphs of the theorem.
We first prove the inequality
$$ (a_{1} + \dots a_{n})^{r} \leq 2^{rn} (a_{1}^{r} + \dots a_{n}^{r} )
\tag 6.9$$
where $a_{1}, \dots a_{n} \geq 0$ and $n \in \N$, and $r$ is some fixed
natural number.
 Recall from the argument after \cite{22}, (3.2.5) that
$(a + b )^{r} \leq 2^{r}(a^{r} + b^{r})$.
Hence, assuming (6.9) holds for $n-1$,
$$ \split(a_{1} + \dots a_{n} &)^{r}
\leq 2^{r}\bigl( (a_{1} + \dots a_{n_1})^{r}
+ a_{n}^{r} \bigr)\\
&
\leq 2^{r}2^{r(n-1)} (a_{1}^{r} + \dots a_{n}^{r}),\endsplit
$$
which proves (6.9) by induction on $n$.
\par
 For convenience, we replace $\tau$ with the
equivalent subadditive scale $1 + \tau$,
so that we have $\tau \geq 1$ and $\tau$ subadditive.
\par
Assume that $A$ is strongly spectral invariant in $B$.
We verify the strong spectral invariance of
$L_{1}^{\tau}(G, A)$ in $L_{1}(G, B)$.  Let
$\bigl\{ \pa \quad \pa_{m} \bigr\}$
be a family of
increasing seminorms for $A$.
We topologize $L_{1}^{\tau}(G, A)$ by the
family of increasing seminorms
$$\parallel \psi \parallel^{\prime}_{m}
=\parallel \psi \parallel_{m,m}.\tag 6.10$$
(See (6.2).)
Note that $\pa \quad \pa_{0}$ is the norm on $L_{1}(G, B)$.
We show that these seminorms satisfy (1.3).
\par
Let $\psi_{1}, \dots \psi_{n} \in G\rtimes^{\tau}A$.
To prepare to
estimate  $\parallel \psi_{1}*\dots \psi_{n} \parallel_{m, m}$,
we write $\psi_{1} *\dots \psi_{n} (g)$ as
$$ \int \dots \int \alpha_{\eta_{1}}(\psi_{1}(h_{1}))
\dots \alpha_{\eta_{n-1}}(\psi_{n-1}(h_{n-1}))
\alpha_{\eta_{n}}(\psi_{n}(h_{n})) dh_{1}\dots dh_{n-1} \tag 6.11 $$
where $h_{1}, \dots, h_{n-1}$ are the variables of integration,
$h_{n} = h^{-1}_{n-1} \dots h_{1}^{-1} g$, $\eta_{1}= e$,
$\eta_{k}= h_{1}\dots h_{k-1}$.
We proceed to estimate.  Using (6.11) and the left invariance of
Haar measure,
$$\aligned
& \parallel \psi_{1}*\dots \psi_{n} \parallel_{m, m} =
\int_{G} \tau^{m}(g)
\parallel \psi_{1}*\dots \psi_{n}(g) \parallel_{m}dg\\
&\leq
 \int \dots \int \tau^{m}(g) \parallel \alpha_{\eta_{1}}(\psi_{1}(h_{1}))
\\ & \qquad
\dots         \alpha_{\eta_{n-1}} (\psi_{n-1}(h_{n-1}))
\alpha_{\eta_{n}}(
\psi_{n}(h_{n}))\parallel_{m}dh_{1}\dots dh_{n-1}
dh_{n}
\endaligned \tag 6.12 $$
Since $A$ is strongly spectral invariant in $B$, we
may bound the normed expression in the integrand of (6.12):
$$\split
\parallel \alpha_{\eta_{1}}(\psi_{1}(h_{1}))  &
\dots
\alpha_{\eta_{n}}(\psi_{n}(h_{n}))\parallel_{m}\leq \\
&D_{1}{C_{1}}^{n}
\sum_{k_{1} + \dots k_{n} \leq p }
\biggl\{\parallel \alpha_{\eta_{1}}(\psi_{1}(h_{1})) \parallel_{k_{1}}
\dots
\parallel
\alpha_{\eta_{n}}(\psi_{n}(h_{n}))\parallel_{k_{n}}\biggr\}
\endsplit \tag 6.13$$
for some constants $D_{1}, C_{1}>0$ and $ p\geq m$
all depending only on $m$.
If $k_{i} \not= 0$, then the temperedness of
the action of $G$ on $A$ gives
$$\pa \alpha_{\eta_{i}}(\psi_{i} (h_{i}))\pa_{k_{i}}
\leq D_{2}\tau^{d}(\eta_{i}) \pa \psi_{i}(h_{i}) \pa_{s}
\tag 6.14  $$
where $s$, $D_{2}$ and $d$ depend
only on $p$, since $k_{i} \leq p$.  If $\{ i_{1}, \dots
i_{p} \}$  contains all the $i_{j}'s$ for which
$k_{i_{j}} \not= 0$ (note that
there are at most $p$ because $k_{1} +
\dots k_{n} \leq p$) in the bracketed
expression of (6.13), then that expression
is bounded by
$$ D_{2}^{p}\tau^{d}(\eta_{i_{1}})
\dots \tau^{d}(\eta_{i_{p}})
\pa \psi_{1}(h_{1})\pa_{k_{1}} \dots \pa \psi_{i_{1}} (h_{i_{1}}) \pa_{s}
\dots \pa \psi_{i_{p}}(h_{i_{p}}) \pa_{s} \dots \pa \psi_{n}(h_{n})
\pa_{k_{n}}, \tag 6.15 $$
where we used the temperedness condition (6.14) $p$ times,
and the fact that $\alpha$ leaves the norm $\pa \quad \pa_{0}$ on
$B$ invariant.
By the subadditivity of $\tau$, we have
$$\split\tau^{d}(\eta_{k}) &\leq (\tau(h_{1}) + \dots \tau(h_{k-1}))^{d}\\
& \leq (\tau(h_{1}) + \dots \tau(h_{n}))^{d}. \endsplit \tag 6.16$$
So (6.13) is bounded by
$$ D_{1}C_{1}^{n}D_{2}^{p} (\tau(h_{1}) + \dots \tau(h_{n}) )^{dp}
\sum_{k_{1} + \dots k_{n} \leq p + ps }
\biggl\{
\pa \psi_{1}(h_{1}) \pa_{k_{1}} \dots \pa \psi_{n}(h_{n})\pa_{k_{n}}
\biggr\}.
\tag 6.17 $$
Plugging this bound on  (6.13) back
into the integrand of (6.12), and also using (6.16) with
$\eta_{i}$ replaced with $g$ and $d=m$, we see
that the integrand of (6.12) is bounded by
$$ \split D_{1} D_{2}^{p} &C_{1}^{n}
(\tau(h_{1}) + \dots \tau(h_{n}))^{dp  + m}
\\
& \sum_{k_{1} + \dots k_{n}\leq p+ps} \biggl\{
\pa \psi_{1}(h_{1}) \pa_{k_{1}} \dots \pa \psi_{n}(h_{n})\pa_{k_{n}}
\biggr\}\endsplit \tag 6.18 $$
which by (6.9) is bounded by (here $C_{3} = 2^{r}$ from (6.9)
with $r= dp  +m$)
$$ \aligned&  D_{1} D_{2}^{p} (C_{1}C_{3})^{n}
\bigl(\tau(h_{1})^{dp+m} + \dots \tau(h_{n})^{dp+m}\bigr)
\\
&\qquad \qquad \qquad \qquad \times
\sum_{k_{1} + \dots k_{n}\leq p+ ps} \biggl\{
\pa \psi_{1}(h_{1}) \pa_{k_{1}} \dots \pa \psi_{n}(h_{n})\pa_{k_{n}}
\biggr\}\\
&\leq  D_{4} C_{4}^{n}
\sum_{k_{1} + \dots k_{n}\leq t}\biggl\{
\tau(h_{1})^{k_{1}}\pa \psi_{1}(h_{1}) \pa_{k_{1}}
\dots \tau(h_{n})^{k_{n}}\pa \psi_{n}(h_{n})\pa_{k_{n}}
\biggr\}\endaligned \tag 6.19 $$
where $t= p + ps + dp  + m$.  Therefore (6.12) is bounded by
$$\aligned \pa& \psi_{1} * \dots \psi_{n} \pa_{m, m}\\
& \leq D_{4} C_{4}^{n} \sum_{k_{1} + \dots k_{n}\leq t}
\int\dots \int
\tau^{k_{1}}(h_{1})\pa \psi_{1}(h_{1}) \pa_{k_{1}}
\dots \tau^{k_{n}}(h_{n})\pa\psi_{n}(h_{n})\pa_{k_{n}}dh_{1}\dots dh_{n}\\
& = D_{4}C_{4}^{n}\sum_{k_{1} + \dots k_{n}\leq t}
\pa \psi_{1}\pa_{k_{1}, k_{1}} \dots \pa \psi_{n} \pa_{k_{n}, k_{n}}
\\
& = D_{4}C_{4}^{n}\sum_{k_{1} + \dots k_{n}\leq t}
\pa \psi_{1}\pa^{\prime}_{k_{1}} \dots \pa \psi_{n} \pa^{\prime}_{k_{n}},
\endaligned \tag 6.20 $$
which, since  $D_{4}$ and $C_{4}$ do not depend on $n$, gives the
strong spectral invariance of $L_{1}^{\tau}(G, A)$ in $L_{1}(G, B)$.
\vskip\baselineskip
\par
The strong spectral invariance of $G\rtimes^{\tau}A$ in $L_{1}(G, B)$
follows from the estimate (6.20) and the following lemma.
\proclaim{Lemma 6.21}
Let $\tau$ be any sub-polynomial scale on $G$.
Topologize $G\rtimes^{\tau}A$ by the increasing seminorms
$$\parallel \psi \parallel^{\prime\prime}_{m}
=\sum_{|\gamma|\leq m}
\parallel \psi \parallel_{m,\gamma,m}.\tag 6.22$$
If $G$ acts differentiably on $A$, then $G\rtimes^{\tau}A$
is a dense right ideal in $L_{1}^{\tau}(G, A)$ and
moreover for all $m \in \N$ there exists
$D>0$ and $k,l\in \N$ such that
$$\pa \varphi * \psi \pa^{\prime\prime}_{m}
\leq D \pa \varphi \pa^{\prime\prime}_{k} \pa
\psi \pa^{\prime}_{l},\tag 6.23$$
for all $\varphi, \psi \in G\rtimes^{\tau}A$.  Similarly, if
 $\tau_{-}$ bounds $Ad$ on $G$, then $G\rtimes^{\tau}A$
is a dense left ideal in $L_{1}^{\tau}(G, A)$, and we have
the inequalities
(6.23) holding, but with $^{\prime\prime}$
and $^{\prime}$ switched on the
right hand side.
\endproclaim
\demo{Proof} The second statement
is just \cite{22}, (2.2.7).
Assume that $G$ acts differentiably on $A$, and replace $\tau$ with
an equivalent scale satisfying $\tau \geq 1$.  We have
$$\aligned &\pa X^{\gamma} \varphi * \psi (g)\pa_{m}
=\pa X^{\gamma} \int_{G} \varphi(gh) \alpha_{gh}(\psi(h^{-1})) dh\pa_{m}\\
&\leq K\sum_{\beta + {\tilde \beta}= \gamma}
\int_{G}\pa X^{\beta}\varphi(gh) \pa_{p}\,\, \pa (X^{\tilde \beta}{\text{
on $g$ }}) \alpha_{gh} (\psi(h^{-1})) \pa_{q} dh
\qquad{\text{prod rule, $A$ Fr\'ech alg}} \\
&\leq {\tilde K}\sum_{\beta + {\tilde \beta}= \gamma}
\int_{G}\pa\tau^{d} X^{\beta}\varphi(gh) \pa_{p}
\,\,
\pa (\psi(h^{-1})) \pa_{r} dh.
\qquad{\text{action diff, then tempered}} \endaligned
\tag 6.24 $$
To see the inequality (6.23), place a $\tau^{m}(g)$ in front of (6.24),
use that $\tau$ is sub-polynomial, and integrate over $g \in G$.
\qed\enddemo
This proves Theorem 6.7.
\qed
\enddemo
\subheading{Remark 6.25}  The strong spectral invariance of
 $G\rtimes^{\tau}B^{\infty}$
in $L_{1}(G, B)$   generalizes Bost \cite{4}, Theorem 2.3.3(a), which
proves a similar theorem for  the case
of elementary Abelian groups.
See also the last example in \cite{22}, \S 5.
\par
For a specific example, let $G$ be $SL_{2}(\Z)$.  Then $G$ is discrete,
finitely generated and Type R,
but does not have polynomial growth.  Let $\tau $ be the word gauge.
The group $G$ has a natural action on the irrational rotation algebra
$A_{\theta}$ \cite{22}, end of \S 5,
and Theorem 6.7  tells us that
$G \rtimes^{\tau} A_{\theta}= L_{1}^{\tau}(G, A_{\theta})$
is strongly spectral invariant in
$L_{1}(G, A_{\theta})$.
\subheading{Example 6.26}
Next let $H$ be a compactly generated polynomial
growth Type R Lie group
(for example a closed subgroup of a connected
nilpotent Lie group).  Let $G$ and $K$ be closed subgroups
of $H$, and let $\tau$ be the word gauge on $H$.
Then $\tau$ restricts to  a gauge
on $G$ which bounds $Ad$ on $H$ \cite{22}, Corollary 1.5.12.
Define a scale $\s$ on $H/K$ as in
Example 2.6.  Let $M = H/K$.  Then $(M, \s, H)$
is a scaled $(G, \tau)$-space.
We may form the nuclear Fr\'echet algebra $G\rtimes^{\tau} \schwmh$,
which is strongly spectral invariant in $L_{1}(G, C_{0}(M))$ by
Theorem 6.7.
\subheading{Example 6.27}
Alternatively, in the preceding example
$(M, \s, G)$ is a scaled $(G, \tau)$-space, and
the (possibly nonnuclear) Fr\'echet algebra
$G\rtimes^{\tau} \Cal S_{G}^{\s}(M)$
is strongly spectral invariant in $L_{1}(G, C_{0}(M))$.
\par
To give a familiar example, let $H=G= \Z$, $K= \{ e\}$.  Then
$\tau $ is the absolute value function on $\Z$ and $G\rtimes^{\tau}
\schwmh$ is $\Z \rtimes \Cal S(\Z)$ with $\Z$ acting by translation,
which is isomorphic to the smooth compact operators defined in \S 5.
\par
For other examples, see \cite{22}, \S 5, or
Example 7.20 below.  We remark that
Theorem 6.7 gives another proof that
strong spectral invariance is preserved
by taking crossed products
with finite groups (see Corollary 4.2(1)).
\vskip\baselineskip
\subheading{Remark 6.28}
In the definition of
$G\rtimes^{\tau} A$ (6.4), we could have let the differential operator
$X^{\gamma}$ act via right translation $g \varphi(h) = \varphi(hg)$
instead of left translation.
Then $G\rtimes^{\tau}A$ would be a Fr\'echet algebra as long as
$\tau$ is sub-polynomial, and the action of $G$ is $\tau$-tempered,
with {\it no} requirement about $\tau_{-}$ bounding $Ad$ or
the action of $G$ being differentiable.  The proof of (the appropriate
modified versions of \cite{22},
Theorem 2.1.5 and Theorem 2.2.6 and) Theorem 6.7
would still go through to give a strongly spectral invariant
smooth crossed product
$G\rtimes^{\tau}A$ in $L_{1}(G, B)$ if $A$ is strongly spectral
invariant in $B$ (in fact, the only changes
in Theorem 6.7  would be in the
estimate (6.24) of Lemma 6.21).  This
would allow us to form strongly
spectral invariant
smooth
crossed products in many cases when
$G$ has no gauge which bounds $Ad$ (for
example, when $G$ is not Type R).
\par
One shortcoming of this approach is that the left \lq\lq covariant
differentiable representations\rq\rq\ of $(G, A)$
would not necessarily be in one to one correspondence with
left \lq\lq differentiable representations\rq\rq\  of the
crossed product $G\rtimes^{\tau}A$, since the action of $G$ on
$G\rtimes^{\tau}A$ on the
left would not necessarily be differentiable (see
\cite{23}, \S 5 and Theorem 5.3).
Another shortcoming is that if $A$ is a *-algebra, in order for the
crossed product $G\rtimes^{\tau}A$
to be a *-algebra, we must require $\tau_{-} \thicksim \tau$,
$\tau_{-}$ bound $Ad$,
and that the action of $G$ on $A$ be
differentiable \cite{22}, \S 4 and
Corollary 4.9.
So,
if we want *-algebras, we may as well stick with
our original definition of $G\rtimes^{\tau}A$.
\par
We remark that if $\tau_{-}$ ($\tau$) bounds $Ad$ on $G$, then
left (right) differentiable operators can be turned into right (left)
ones.
\vskip\baselineskip
\par
We investigate what happens when $G$ is not Type R.  We already know by
\cite{22}, \S 1.4
that $G$ has no gauge that bounds $Ad$.
In fact, we have the following more decisive result.

\proclaim{Theorem 6.29} Let $G$ be any Lie group,
and let $\tau$ be any sub-polynomial
scale on $G$. Assume also that
$\tau_{-} \thicksim \tau$ and that $\tau_{-}$
bounds $Ad$.  Under these conditions, if the
Lie group $G$ is {\it not} Type R, then the group algebra
$\Cal S^{\tau}(G)$ is {\it never} spectral invariant
in $L_{1}(G)$, or in either of the C*-algebras
$C_{r}^{*}(G)$ or $C^{*}(G)$, whatever the choice of
$\tau $ satisfying the above conditions.
\endproclaim
\par
The theorem basically says that,
via all the ways of showing that
$S^{\tau}(G)$ is a Fr\'echet *-algebra
I know of, in general we never get spectral invariance if $G$ is not Type R.
\demo{Proof}
By \cite{22}, Theorem 5.17, and since $\tau_{-} \thicksim \tau$,
we know that $\Cal S^{\tau}(G)$ is a Fr\'echet
*-algebra.  We define a simple $\Cal S^{\tau}(G)$-module
which is not contained in an $L_{1}(G)$-module.
\par
We define a $G$-module
$V$ by taking $V$ to be the Lie algebra of $G$, and letting
$g v \equiv Ad_{g}v$.
Let $X\in V$ have eigenvalue $\lambda$ not on the unit circle for
some $g \in G$.  If $V$ has a non-trivial invariant subspace $W$,
then either $X\in W$, or $X$ has nonzero image in $V/W$.
Continuing, we eventually reach a simple, finite
dimensional $G$-module $\tilde V$ with an eigenvector $\tilde X$ with
eigenvalue $\lambda$.
\par
Since $\tau_{-} \thicksim \tau$, we know $\tau$ bounds $Ad$.
Thus we may
integrate the original representation of $G$ on $V$ to a representation
of $A=\Cal S_{1}^{\tau}(G)$. This representation of $A$
goes through the argument of the previous paragraph to give an
irreducible representation of
$A$ on $\tilde V$.
\par
By \cite{21}, Theorem 1.4, if $A$ is spectral invariant
in any of the  algebras $L_{1}(G)$, $C_{r}^{*}(G)$, or
$C^{*}(G)$ - call them $B$ - then
$\tilde V$ must have a continuous irreducible
$B$-module structure extending the
action of $A$ on $\tilde V$.
\par
Clearly for $C^{*}(G)$ this is impossible since this would imply that
the original representation of $G$ on $\tilde V$
were unitary, contradicting
that $g$ has eigenvalue $\lambda$.  We show that in fact $\tilde V$
cannot have an $L_{1}(G)$-module structure (which proves the theorem).
For assume that it does.  Then $G$ acts as left multipliers on $L_{1}(G)$
via $k\varphi(h) = \varphi(k^{-1}h)$, and we have $\parallel
 k \parallel_{mult}= 1$
for all $k \in G$.  If $w$ is a nonzero element of
$ \tilde V$, any $v \in \tilde V$
can be written $\varphi w$
for some $\varphi \in L_{1}(G)$.  We then define an action
of $G$ on $\tilde V$ by $k v = (k\varphi)w$, which must agree with
the original representation of $G$ on $V$.  Since
$\parallel k \parallel_{mult} =1$,
there is some constant
$C$ such that $\parallel k \parallel_{\Cal B({\tilde V})}\leq C$.
But for large $n$, $g^{n}$ has arbitrarily
large eigenvalues, and hence arbitrarily large
norm in $\Cal B({\tilde V})$.
This is a contradiction, so there can be no $L_{1}(G)$-module
structure on $\tilde V$.
\qed \enddemo
\heading \S 7 Crossed Products by Polynomial Growth Groups
\endheading
\par
We use the results of \S 6 to show that our dense subalgebras
are spectral invariant in the C*-crossed
product $G\rtimes B$, at least when
$G$ has polynomial growth.
\par
Throughout this section, unless otherwise stated,
$B$ will be a C*-algebra.
We shall use $\pa \quad \pa$ to denote the norm on the reduced
C*-crossed
product $G\rtimes_{r} B$.
By Paterson \cite{16}, Proposition 0.13,
$G$ is amenable if $G$ has polynomial growth
(see Definition 2.5).
So in this case, $G\rtimes B = G\rtimes_{r} B$.
 If $G$ is compactly generated,
let $\tau $ be the word gauge on $G$ (see Definition 2.5).
We replace  $\tau$ with $1 + \tau$, so $\tau  \geq 1$
and $\tau $ is subadditive and submultiplicative (namely
$\tau(gh)\leq \tau(g) + \tau(h)$ and $\tau(gh) \leq \tau(g) \tau(h)$).
Let $L_{1}(G, \tau^{q})$ and $L_{1}(G, B, \tau^{q})$
be the Banach *-algebras of $L_{1}$ functions corresponding to the
measure $\tau^{q} dg$ on G.
Note that by Theorem 6.7 (or the estimates in its proof),
we know that $L_{1}(G, B, \tau^{q})$ is
strongly spectral invariant in $L_{1}(G, B)$, so we have equality
of the two spectral radii
$$ \lim_{n\rightarrow \infty} \pa \varphi^{n} \tau^{q} \pa_{1}^{1/n}
=  \lim_{n\rightarrow \infty} \pa \varphi^{n}
\pa_{1}^{1/n}= \nu(\varphi), \qquad
\varphi \in L_{1}(G, B, \tau^{q}).\tag 7.1 $$
Here we let $\nu(\varphi)$ denote the spectral radius of $\varphi$ in
$L_{1}(G, B)$.
We shall imitate the argument of Pytlik \cite{19},
generalizing the proof from
the case $L_{1}(G, \tau^{q})$ to the case $L_{1}(G, B, \tau^{q})$,
and then use this to show that $L_{1}(G, B, \tau^{q})$ is
spectral invariant in the C*-crossed product $G\rtimes B$.
We shall be using the star operation
$f^{*} (g) = \Delta(g) \alpha_{g}(f(g^{-1})^{*})$ for functions
$f \colon G \rightarrow B$.
\proclaim{Lemma 7.2 (Compare Pytlik \cite{19}, Lemma 4)}
Let $G$ be any locally compact group.  Let $\varphi =
\varphi^{*} \in L_{1}(G, B)$, and let $D$ be any dense subset of
$L_{1}(G, B)$.  Then there exists   functions $ f_{1}^{*}, f_{2} \in D$
such that
$$ \nu(\varphi ) \leq \limsup_{n \rightarrow \infty}
\pa f_{1}*\varphi^{n} * f_{2} \pa_{1}^{1/n}. \tag 7.3$$
\endproclaim
\demo{Proof}
For $\varphi \equiv 0$, this is clear.  Let $ \varphi \not= 0$.
Then $\nu (\varphi) \not= 0$, since, for example, $L_{1}(G, B)$
has a faithful *-representation on a Hilbert space.
\par
Let $a_{n} = \pa \varphi^{n+2} \pa_{1} \pa \varphi^{n}\pa_{1}^{-1}$.
Then $\limsup_{n \rightarrow \infty} a_{n} =
\lim_{n\rightarrow \infty} \pa
\varphi^{n}\pa_{1}^{2/n} = \nu(\varphi)^{2}\not= 0$
by  manipulations of limits of positive real numbers.
\par
Let $0< a< \limsup_{n\rightarrow \infty} a_{n}$.
Let $\epsilon < \max (1, a (6\pa \varphi \pa_{1}
)^{-1}) $, and choose
$f_{1} \in D^{*} $, $f_{2} \in D$
such that $\pa f_{i} - \varphi \pa_{1} < \epsilon$.
Then $\pa f_{i} \pa_{1} \leq	 2 \pa \varphi \pa_{1}$ and we have
$$
\aligned \pa \varphi^{n+2} \pa_{1} & \leq \epsilon
\pa \varphi^{n+1} \pa_{1}
+ \pa f_{1} * \varphi^{n+1} \pa_{1} \\
& \leq \epsilon \pa \varphi^{n+1}
\pa_{1} + \epsilon \pa f_{1}* \varphi^{n} \pa_{1}
+ \pa f_{1} * \varphi^{n} * f_{2} \pa_{1}  \\
& \leq
a/2\pa \varphi^{n} \pa_{1} + \pa f_{1} * \varphi^{n} * f_{2} \pa_{1}.
  \endaligned
$$
Hence
$$ \pa f_{1} * \varphi^{n} * f_{2} \pa_{1}   \geq \pa \varphi^{n+2} \pa_{1}
- a/2 \pa \varphi^{n}
\pa_{1} = \pa \varphi^{n} \pa_{1} \bigl( a_{n}-a/2\bigr). $$
Since $a_{n} - a/2 \geq a/2$ for infinitely many $n$,
$$ \limsup_{n\rightarrow \infty} \pa f_{1} * \varphi^{n} * f_{2}
 \pa_{1}^{1/n} \geq
\lim_{n \rightarrow \infty}
\pa \varphi^{n} \pa_{1}^{1/n} = \nu (\varphi).$$
\qed
\enddemo
\proclaim{Lemma 7.4 (Compare Bost \cite{4}, (7.3.10))}
Let $G$ be any unimodular locally compact group.
Let $D$ be the
vector space of measurable, compactly supported,  step functions from
$G$ to $B$.  Then there exists a  norm $\pa \quad \pa_{D}$
on $D$ such that
$$ \pa f_{1} * \psi * f_{2}
\pa_{\infty} \leq \pa f^{*}_{1} \pa_{D}\, \pa \psi \pa\,
\pa f_{2} \pa_{D},
\qquad f^{*}_{1}, f_{2} \in D, \quad \psi \in L_{1}(G, B). \tag 7.5$$
\endproclaim
\demo{Proof}
We first consider the case $ f^{*}_{1} = \xi_{1} \otimes b^{*}_{1}$,
$f_{2} = \xi_{2} \otimes b_{2}$, where
$\xi_{i} $ is a characteristic function of a relatively
compact measurable subset
of $G$,  and $b_{i} \in B$.
For $g \in G$, we bound $\pa (f_{1} * \psi * f_{2})(g) \pa_{B}$.
Let $B$ be faithfully *-represented on a Hilbert space $\H$, and
let $\eta_{1}, \eta_{2} \in \H$.  By changes of variables,  we have
$$  \split < \eta_{1}, (f_{1} * \psi * f_{2} )(g) \eta_{2} >_{\H}
= <\eta_{1}, (\xi_{1}^{*} & *(b_{1} \psi b_{2} ) *
\xi_{2})(g) \eta_{2} >_{\H}
\\&= < \xi_{1} \otimes \eta_{1} , (b_{1} \psi b_{2})((\xi_{2})_{g} \otimes
\eta_{2}) >_{L_{2}(G, \H )}, \endsplit $$
where $(\xi_{2})_{g} \in L_{2}(G)$ is the function $(\xi_{2})_{g} (h)
= \xi_{2}(hg)$, and $(b_{1} \psi b_{2}) \in L_{1}(G, \H )$
acts on $((\xi_{2})_{g} \otimes \eta_{2} ) \in L_{2}(G, \H )$
via the regular representation induced from the representation of
$B$ on $\H$ Pedersen \cite{17}, \S 7.7.
It follows that
$$\aligned |<\eta_{1}, (f_{1} * \psi * &f_{2} )(g) \eta_{2} >_{\H}|\leq
\pa \xi_{1} \otimes \eta_{1} \pa_{L_{2}(G, \H)}\,
\pa b_{1} \psi b_{2} \pa \,\pa (\xi_{2})_{g} \otimes
\eta_{2} \pa_{L_{2}(G, \H)}\\
&
= \biggl(\pa \xi_{1} \pa_{L_{2}(G)}\, \pa b_{1} \psi b_{2} \pa\,
\pa \xi_{2} \pa_{L_{2}(G)} \biggr)\,\pa \eta_{1} \pa_{\H}\,
\pa \eta_{2} \pa_{\H} \\
& \leq
\biggl(\pa \xi_{1} \pa_{L_{2}(G)}\,\pa b_{1} \pa_{B} \,\pa  \psi  \pa\,
\pa \xi_{2} \pa_{L_{2}(G)} \,\pa b_{2} \pa_{B}\biggr)\,
\pa \eta_{1} \pa_{\H} \,\pa \eta_{2} \pa_{\H}. \endaligned \tag 7.6$$
Let $\pa \quad \pa_{D}$ be the norm on $D$ inherited from the
projective tensor product $ L_{2}(G) \tensp B$.
Then by (7.6) and the definition of the projective topology, we have
$$ |<\eta_{1},  (f_{1} * \psi * f_{2} )(g) \eta_{2} >_{\H}| \leq
\biggl(\pa f_{1}^{*} \pa_{D}\, \pa  \psi  \pa\,
\pa f_{2}\pa_{D}\biggr)\, \pa \eta_{1} \pa_{\H} \,\pa \eta_{2} \pa_{\H},
\tag 7.7 $$
for  $f^{*}_{1}$ and $f_{2}$ in $ L_{2}(G) \tensp B$.
Taking the sup over $\pa\eta_{i}\pa_{\H} \leq 1$
we have
$$ \pa (f_{1} * \psi * f_{2} )(g) \pa_{B} \leq
\pa f_{1}^{*}\pa_{D} \, \pa  \psi  \pa\,
\pa f_{2} \pa_{D}.\tag 7.8
$$
Taking the sup over $g \in G$, we get (7.5).
\qed
\enddemo
\proclaim{Lemma 7.9 (Compare Pytlik \cite{19}, Lemma 5)}
Let $G$ be a compactly generated polynomial growth
group.  Let $D$ be as in Lemma 7.4, and let $f_{1}^{*} \in D$,
$f_{2} \in D$, and $ \psi \in L_{1}(G, B, \tau^{q})$.
Then there exists constants $M, N>0$ (not depending on $\psi$)
such that for $m \in \N^{+}$
$$ \pa f_{1} * \psi * f_{2} \pa_{1} \leq \pa \psi \pa M m^{r} +
\pa  \psi\tau^{q} \pa_{1} Nm^{-q},\tag 7.10$$
where $r>0$ is the growth constant of the group.
\endproclaim
\demo{Proof} Let $U$ be a generating set for $G$.
We have $\pa f_{1} * \psi * f_{2} \pa_{1} = \pa (f_{1}* \psi * f_{2})
\chi_{U^{m}} \pa_{1}
+ \pa (f_{1} * \psi * f_{2} ) \chi_{G - U^{m}} \pa_{1}$,
where $\chi_{U^{m}}$ and $\chi_{G- U^{m}}$ are characteristic
functions of the sets $U^{m}$ and $G - U^{m}$ respectively.
But by Lemma 7.4 and
since $G$ is unimodular \cite{16}, Proposition 6.9,6.6,
$$ \pa (f_{1} * \psi * f_{2}) \chi_{U^{m}} \pa_{1}
\leq \pa f_{1}* \psi * f_{2} \pa_{\infty}\, \pa \chi_{U^{m}} \pa_{1}\,
\leq \pa f^{*}_{1} \pa_{D}\, \pa f_{2}\pa_{D}\,
\pa \psi \pa m^{r} = M \pa \psi \pa m^{r}. $$
Also
$$ \split \pa (f_{1} * \psi * f_{2} ) \chi_{G- U^{m}} \pa_{1}
\leq &\pa (f_{1} * \psi * f_{2}) \tau^{q} \pa_{1}
\pa \chi_{G- U^{m}} \tau^{-q} \pa_{\infty}  \\
&
\leq	 \pa \psi \tau^{q} \pa_{1}\, \pa f_{1} \tau^{q} \pa_{1}
\pa f_{2} \tau^{q} \pa_{1}
(1+m)^{-q} \leq \pa \psi \tau^{q} \pa_{1} N m^{-q}. \endsplit $$
\qed
\enddemo
\proclaim{Theorem 7.11 (Compare Pytlik \cite{19}, Theorem 6)}
Let $G$ be a compactly generated polynomial growth group.
For every $\varphi = \varphi^{*} \in L_{1}(G, B, \tau^{q})$
we have $\nu(\varphi) = \pa \varphi \pa $. \endproclaim
\demo{Proof}
It suffices to prove $\nu(\varphi) \leq \pa \varphi \pa$.
Let $a \geq 1$ be arbitrary, and let $m_{n}$ be a sequence
of integers such that $\lim_{n \rightarrow \infty} m_{n}^{1/n} =a$.
Putting $\varphi^{n}$ for $\psi$ and $m_{n}$ instead of $m$ in (7.10)
we get
$$ \pa f_{1} * \varphi^{n} * f_{2} \pa_{1}^{1/n} \leq
\biggl( \pa \varphi \pa^{n} Mm_{n}^{r} + \pa \varphi^{n} \tau^{q} \pa_{1}
Nm_{n}^{-q} \biggr)^{1/n}.\tag 7.12$$
If $n$ tends to infinity, the right side of (7.12) tends to a limit,
which by (7.1) is equal to $\max\bigl\{\pa\varphi\pa a^{r}, \nu(\varphi)
a^{-q} \bigr\}$.  Therefore
$$ \limsup_{n\rightarrow\infty} \pa f_{1} * \varphi^{n} *f_{2} \pa_{1}^{1/n}
\leq \max\biggl\{ \pa \varphi \pa a^{r}, \nu(\varphi)
a^{-q} \biggr\},$$
which for $a = (\pa \varphi \pa^{-1} \nu(\varphi))^{1\over{r+q}} \geq 1$
and for $f_{1}$ and $f_{2}$ as in Lemma 7.4  yields
$$\nu(\varphi) \leq \pa \varphi \pa^{q\over{r+q}}
\nu(\varphi)^{r\over{r+q}}$$
and so
$$\nu(\varphi) \leq \pa \varphi \pa.$$
This proves Theorem 7.11
\qed
\enddemo
\par
The following theorem is essentially Hulanicki \cite{8}, Proposition 2.5.
\proclaim{Theorem 7.13}
Let $\Cal A$ be a Banach *-algebra.
Assume that $\Cal A$ is faithfully *-represented in $\Cal B(\H)$ for
some Hilbert space $\H$, such that the C*-norm $\pa a \pa_{\Cal B(\H)}$
is equal to the spectral radius of $a$ in $\Cal A$, for all $a= a^{*}$
in $\Cal A$.  Then
$\text{spec}_{\Cal A}(a) = \text{spec}_{\Cal B(\H)} (a) $for
all $a \in \Cal A$.
\endproclaim
\demo{Proof}  Let $\Cal B$ be the closure of $\Cal A$ in $\Cal B(\H)$.
We wish to show that $\Cal A$ is spectral invariant in $\Cal B$.
By Hulanicki \cite{8}, Proposition 2.5, we have that
$\text{spec}_{\Cal A}(a) = \text{spec}_{\Cal B} (a)$ for
all $a=a^{*}$  in  $\Cal A$.
It follows that for every $a= a^{*}$ in
$\tilde \Cal A$, $a$ is invertible in
$\tilde \Cal A$ iff $a$ is invertible in $\tilde \Cal B$.
\par
Assume for
a contradiction that $\Cal A$ is not
spectral invariant in $\Cal B$.  Then by \cite{21}, Theorem 1.4, there
is a maximal left ideal $I$
in $\tilde\Cal  A$ which is dense in $\tilde\Cal B$.
Hence $I$ contains an invertible element $a$ of $\tilde \Cal B$ which is
not invertible in $\tilde \Cal A$.  But then $a^{*}a$ is a self-adjoint
element of $\tilde \Cal A$
which is in $I$ and hence not invertible in $\tilde\Cal  A$,
but is invertible in $\tilde\Cal  B$.
This is a contradiction and completes the
proof.
\qed
\enddemo
\proclaim{Corollary 7.14} Let $G$ be a compactly generated
polynomial growth group. Then the
Banach *-algebra $L_{1}(G, B, \tau^{q})$ is
spectral invariant in the C*-crossed product $G \rtimes B$
for any $ q \in \N$.
Hence $L_{1}^{\tau}(G, B)$ is spectral invariant in $G \rtimes B$.
\endproclaim
\demo{Proof}
By Theorem 7.11 and (7.1), the spectral radius of a self-adjoint
element of $L_{1}(G, B, \tau^{q})$
is equal to it's norm in $G\rtimes B$.
Hence by Theorem 7.13, $L_{1}(G, B, \tau^{q})$ is spectral invariant
in $G\rtimes B$.  Since $$L_{1}^{\tau}(G, B)= \cap_{q \in \N}
L_{1}(G, B, \tau^{q}),$$
we have the spectral invariance of $L_1^\tau (G, B)$ in $G\rtimes B$.
\qed
\enddemo
\par
The following corollary generalizes Bost \cite{4}, Corollary 2.3.4, which
gives the same result for elementary Abelian groups.
\proclaim{Corollary 7.15}
Let $G$ be a compactly generated polynomial growth group.
 Then the inclusion map $L_{1}(G, B) \hookrightarrow
G\rtimes B$ is an isomorphism of $K$-theory
\endproclaim
\demo{Proof} By  Corollary 7.14, $L_{1}^{\tau}(G, B)$ is
spectral invariant and dense
in $G\rtimes B$.  As we noticed at the beginning of \S 7, it is
also spectral invariant and dense in $L_{1}(G, B)$.
Hence by
\cite{5}, VI.3
and \cite{21}, Lemma 1.2, Corollary 2.3,
or \cite{4}, Appendix,
all three algebras have the same $K$-theory.
\qed \enddemo
\proclaim{Corollary 7.16}
 Let $G$ be a  compactly
generated polynomial growth  group, and
let $\tau$ be the word   gauge on $G$ (Definition 2.5).
Assume that the action of  $G$ on $A$ is
$\tau$-tempered. Then the Fr\'echet algebra $L_{1}^{\tau}(G, A)$
is  spectral invariant
in the C*-crossed product
$G\rtimes B$ if $A$ is strongly spectral invariant in $B$.
\par
If, in addition, $G$ is a Lie group, and either $\tau$ bounds $Ad$
on $G$ or $G$ acts differentiably
on $A$, then $G\rtimes^{\tau} A$ is
a Fr\'echet algebra which is spectral invariant in
$G\rtimes B$ if $A$ is strongly spectral invariant in $B$.
\par
In particular, if $(M, \s, H)$ is a scaled $(G, \tau)$-space, then
the smooth crossed product $G\rtimes^{\tau}\schwmh$ is spectral invariant
in the C*-crossed product $G\rtimes C_{0}(M)$.  Also,
$G\rtimes^{\tau}B^{\infty}$
is spectral invariant in $G\rtimes B$, where $B^{\infty}$
denotes the set of $C^{\infty}$-vectors for the action of $G$ on $B$.
\par
(For these results to apply, $G$ may have to be a Type R Lie
group
- see remarks
preceding  Theorem 6.7, and Theorems 6.7 and 6.29. In general,
a compactly generated polynomial growth Lie
group need not be Type R \cite{11}, Example 1.
However, the word gauge always bounds $Ad$
if such a group is Type R \cite{22}, Corollary 1.5.12.)
\endproclaim
\par
Let $G$ be a compactly generated polynomial growth Lie group.
As noted in the corollary, there are examples when $G$ is not Type R.
By Theorem 6.29 above, we therefore also have examples of $G$-spaces $M$
for which the smooth crossed product $G\rtimes^{\tau}\schwmg$
is never spectral invariant in $G\rtimes C_{0}(M)$, for any
choice of $\s$ and $\tau$ which makes $G\rtimes^{\tau}\schwmg$
a Fr\'echet *-algebra.
However, if $G$ is discrete, a closed subgroup of a connected
polynomial growth Lie group, or if the connected component
of the identity
$G_{0}$ of $G$ is simply connected, then $G$ is Type R \cite{22},
Theorem 1.5.13,
so for large classes of groups we do not have this problem and Corollary
7.16 applies
(see also the examples mentioned in the introduction and the abstract).
\demo{Proof}
 By Theorem 6.7, $G\rtimes^{\tau} A$
and $L_{1}^{\tau}(G, A)$ are strongly spectral invariant
and hence spectral invariant
in $L_{1}(G, B)$.  Also, $L_{1}^{\tau}(G, B)$
is spectral invariant in $L_{1}(G, B)$
by Theorem 6.7, so $G\rtimes^{\tau }A$
and $L_{1}^{\tau}(G, A)$ are both
spectral invariant in $L_{1}^{\tau}(G, B)$.
Since the latter algebra is spectral
invariant in $G\rtimes B$ by  Corollary
7.14, this completes the proof.
\qed
\enddemo
\proclaim{Corollary 7.17}  If $G$ is any compactly
generated polynomial growth Lie group
and $\tau$ is the word gauge on $G$, then the group Schwartz
algebra $\Cal S_{1}^{\tau}(G)$ is spectral invariant
in $C^{*}(G)$.
\endproclaim
\subheading{Remark 7.18}
Corollary 7.17  generalizes the corresponding results Ludwig \cite{12},
Proposition 2.2
for the Schwartz algebra of a nilpotent Lie group, and  Ji \cite{9},
Corollary 1.4
for the Schwartz algebra of a
finitely generated polynomial growth discrete group.
\subheading{Remark 7.19} The statement
about the spectral invariance of $G\rtimes^{\tau} B^{\infty}$
in Corollary 7.16 generalizes the corresponding result
Bost \cite{4}, Theorem 2.3.3(b) for  elementary Abelian groups $G$.
\subheading{Example 7.20}
 Let $H$, $G$ and $K$ be as in Examples 6.26-7, with
$M = H/K$.
Then $G\rtimes^{\tau} \schwmh$ is spectral invariant in the C*-crossed
product $G\rtimes C_{0}(M)$.  Similarly,
$G \rtimes^{\tau} \Cal S^{\s}_{G}(M)$
is spectral invariant in $G\rtimes C_{0}(M)$.
\par
As a special case, the smooth crossed product $\Z \rtimes \Cal S(\Z)$ of
Example 6.27 is spectral invariant in
the compact operators $\Z \rtimes c_{0}(\Z)$.
This also follows from Corollary 5.6 above.
\par
Other examples lie in \cite{22}, \S 5.  For example, if $H$ is
any closed subgroup of $G=GL(n,\R)$ which consists of upper triangular
matrices with $\pm 1$'s on the diagonal, then \cite{22}, Example 5.23 gives
spectral invariant dense subalgebras $H \rtimes^{\tau} \Cal S_{G}^{\s}(M)$
of the C*-crossed product $H \rtimes C_{0}(M)$, where $M= \R^{n}$ and
$H$ and $G$ act by matrix multiplication, or $M = M(n, \R)$ and $H$
and $G$ act by conjugation.
\vskip\baselineskip
\subheading{Remark 7.21} If we defined the smooth crossed product
$G\rtimes^{\tau}A$ with differential operators acting on the right
instead of the left as in Remark 6.28, Corollary 7.16 would still
give the spectral invariance without requiring a gauge that bounds $Ad$
to form the crossed product.
Hence we would have spectral invariant dense subalgebras of
smooth functions for C*-crossed
products by arbitrary compactly generated polynomial growth Lie groups,
with no assumption about $G$ being Type R. However, see the shortcomings
of such algebras mentioned in Remark 6.28.
\subheading{Remark 7.22} We describe an alternate proof of Corollary 7.14
in the case that $G$ is a finitely generated discrete polynomial growth
group.  This proof, along with Pytlik \cite{19},
is what initially suggested to me that Corollary 7.14
might be true.  We show that $L^{\tau}_{1}(G, B)$ is spectral invariant
in $G\rtimes B$. Let $B$ be faithfully *-represented on a Hilbert space
$\H$.  We have the standard representation
$$ \varphi \xi(g) = \sum_{h \in G} \alpha_{g^{-1}}(\varphi(h)) \xi(h^{-1}g),
\qquad
\varphi \in G\rtimes B, \quad \xi \in L_{2}(G, \H). $$
Define a self-adjoint unbounded operator $D$ on $L_{2}(G, \H)$
by
$ D\xi(g) = \tau(g) \xi(g)$.
Define a derivation $\delta$ on $\Cal B(L_{2}(G, \H))$ by
$\delta(T) = i[D, T]$.  Then by Ji \cite{9}, Theorem 1.2, the
set $(G\rtimes B)^{\infty}$ of $C^{\infty}$-vectors for the
action of $\delta$ on $\Cal B(L_{2}(G,\H))$,
which lie in $G\rtimes B$, is a spectral invariant subalgebra
of $G\rtimes B$, as long as it is dense.
It is straightforward to show that
$L_{1}^{\tau}(G, B) \subseteq (G\rtimes B)^{\infty}$, so we have the
density.
We show that $(G\rtimes B)^{\infty} \subseteq L_{1}^{\tau}(G, B)$.
\par
If $\varphi \in (G\rtimes B)^{\infty}$, then we have
$$ \pa (\delta^{k} \varphi)\xi \pa_{L_{2}(G, \H)} \leq
C_{k, \varphi}, $$
for $\pa \xi \pa_{L_{2}(G, \H)}\leq 1$.  Define $\xi=
\delta_{e} \otimes \eta$,
where $\delta_{e}$ is the delta function at $e$, and $\eta \in \H$.  Then
a simple inductive argument shows that
$$ (\delta^{k} \varphi) \xi (g)
= i^{k}\tau^{k}(g) \alpha_{g^{-1}}(\varphi(g) )\eta. $$
This is similar to the formula for $(\delta^{k} \varphi)\xi$ in
\cite{9}, \S 1.
By definition of the $L_{2}$-norm on $L_{2}(G, \H)$, we have
$$C_{k, \varphi}^{2} \geq \pa (\delta^{k}\varphi)\xi \pa_{L_{2}(G, \H)}^{2}
= \sum_{g \in G}\tau^{2k}(g) \pa \alpha_{g^{-1}}(\varphi(g))\eta \pa_{\H}^{2}
\geq \tau^{2k}(g)\pa \alpha_{g^{-1}}(\varphi(g))\eta \pa_{\H}^{2}
$$
for each $g \in G$.
(The last inequality is what uses the discreteness of $G$.)
Taking the sup over $\pa \eta\pa_{\H}\leq 1$, and using the fact that
$\alpha$ is an isometry on $B$, we have
$$\tau^{k}(g)\pa \varphi(g) \pa_{B} \leq C_{k, \varphi}.$$
It follows that for $p \in \N$,
$$ \tau^{k}(g) \pa \varphi(g) \pa_{B}
\leq {C_{k, \varphi} + C_{k +p, \varphi}
\over{1+\tau^{p}(g)}}.$$
Since the right hand side is summable over $g \in G$
for some $p \in \N$ \cite{22}, Proposition 1.5.1,
we have $\varphi \in L^{\tau}_{1}(G, B)$.
So $(G\rtimes B)^{\infty} \subseteq
L_{1}^{\tau}(G, B)$, and the two sets are equal.
By our remarks above,
it follows  that $L_{1}^{\tau}(G, B)$
is spectral invariant in $G\rtimes B$.
\enddocument